\documentclass{amsart}
\usepackage{eepic}
\allowdisplaybreaks
\theoremstyle{plain}
  \newtheorem{theorem}{Theorem}[section]
  \newtheorem{lemma}[theorem]{Lemma}
  \newtheorem{corollary}[theorem]{Corollary}
  \newtheorem{proposition}[theorem]{Proposition}
  \newtheorem{conjecture}[theorem]{Conjecture}
\theoremstyle{definition}
  \newtheorem{definition}[theorem]{Definition}
  \newtheorem{problem}[theorem]{Problem}
\theoremstyle{remark}
  \newtheorem{remark}[theorem]{Remark}
  \newtheorem*{acknowledgements}{Acknowledgements}
\numberwithin{equation}{section}
\newcommand{\CC}{{\mathbb C}}

\newcommand{\ZZ}{{\mathbb Z}}
\newcommand{\RR}{{\mathbb R}}

\newcommand{\A}{{\mathcal A}}
\newcommand{\D}{{\mathcal D}}

\newcommand{\Ell}{{\mathcal L}}

\newcommand{\SplitUnion}
  {\raisebox{-1.00mm}
  {\begin{picture}(4.00,4.00)
  \put(0.00,0.00){\line(1,0){4.00}}
  \put(1.00,0.00){\line(0,1){4.00}}
  \put(3.00,0.00){\line(0,0){4.00}}
  \end{picture}}}

\begin{document}
\title
[A weight system derived from    the               Conway potential function]
{A weight system derived from \\ the multivariable Conway potential function}
\author{Hitoshi Murakami}
\thanks{The author is an EPSRC Visiting Fellow (grant number GR/K/46743)
in the University of Liverpool.\\
{\large {\bf To appear in J. London Math. Soc.}}}
\address{
Department of Mathematics, Osaka City University, 3-138, Sugi\-moto 3-cho\-me,
Sumi\-yoshi-ku, OSAKA 558, JAPAN
\\and\\
Department of Pure Mathematics, University of Liverpool, PO Box 147, LIVERPOOL L69 3BX, UK}
\curraddr{
Department of Mathematics, School of Science and Engineering, Waseda University, 4-1, Oh\-kubo 3-cho\-me, Shin\-juku-ku, TOKYO 169, JAPAN}
\begin{abstract}
A weight system is defined from the (multivariable) Conway potential function.
We also show that it can be calculated recursively by using five axioms.
\end{abstract}
\maketitle
\thicklines
\unitlength 1.00mm

\newcommand{\TwoParallelArcsUpUp}[3]
  {\raisebox{-6.75mm}
  {\begin{picture}( 5.00,15.00)
  \put( 0.00, 0.00){\vector( 0, 1){16.00}}
  \put( 5.00, 0.00){\vector( 0, 1){16.00}}
  {#1}
  \put( 0.00,-0.50){\makebox( 0.00, 0.00)[ct]{#2}}
  \put( 5.00,-0.50){\makebox( 0.00, 0.00)[ct]{#3}}
  \end{picture}}}

\newcommand{\ThreeParallelArcsUpUpUp}[4]
  {\raisebox{-6.75mm}
  {\begin{picture}(10.00,15.00)
  \put( 0.00, 0.00){\vector( 0, 1){16.00}}
  \put( 5.00, 0.00){\vector( 0, 1){16.00}}
  \put(10.00, 0.00){\vector( 0, 1){16.00}}
  {#1}
  \put( 0.00,-0.50){\makebox( 0.00, 0.00)[ct]{#2}}
  \put( 5.00,-0.50){\makebox( 0.00, 0.00)[ct]{#3}}
  \put(10.00,-0.50){\makebox( 0.00, 0.00)[ct]{#4}}
  \end{picture}}}

\newcommand{\ThreeParallelArcsUpUpUpSmall}[4]
  {\raisebox{- 9.50mm}
  {\begin{picture}(10.00,15.00)
  \put( 0.00, 5.00){\vector( 0, 1){11.00}}
  \put( 5.00, 5.00){\vector( 0, 1){11.00}}
  \put(10.00, 5.00){\vector( 0, 1){11.00}}
  {#1}
  \put( 0.00, 4.50){\makebox( 0.00, 0.00)[ct]{#2}}
  \put( 5.00, 4.50){\makebox( 0.00, 0.00)[ct]{#3}}
  \put(10.00, 4.50){\makebox( 0.00, 0.00)[ct]{#4}}
  \end{picture}}}

\newcommand{\OneArcUpAndOneCurlyArcTop}[3]
  {\raisebox{-11.25mm}
  {\begin{picture}(15.00,20.00)
  \put( 0.00, 0.00){\vector( 0, 1){21.00}}
  \put( 5.00, 0.00){\line( 0, 1){12.00}}
  \put(10.00, 5.00){\vector( 0, 1){16.00}}
  \put(10.00,12.00){\oval(10.00, 5.00)[t]}
  \put(12.50, 5.00){\oval( 5.00, 5.00)[b]}
  \put(15.00, 5.00){\line( 0, 1){ 7.00}}
  {#1}
  \put( 0.00,-0.50){\makebox( 0.00, 0.00)[ct]{#2}}
  \put( 5.00,-0.50){\makebox( 0.00, 0.00)[ct]{#3}}
  \end{picture}}}

\newcommand{\OneArcUpAndOneCurlyArcBottom}[3]
  {\raisebox{-11.25mm}
  {\begin{picture}(15.00,20.00)
  \put( 0.00, 0.00){\vector( 0, 1){21.00}}
  \put( 5.00, 5.00){\vector( 0, 1){16.00}}
  \put(10.00, 0.00){\line( 0, 1){12.00}}
  \put(12.50,12.00){\oval( 5.00, 5.00)[t]}
  \put(10.00, 5.00){\oval(10.00, 5.00)[b]}
  \put(15.00, 5.00){\line( 0, 1){ 7.00}}
  {#1}
  \put( 0.00,-0.50){\makebox( 0.00, 0.00)[ct]{#2}}
  \put(10.00,-0.50){\makebox( 0.00, 0.00)[ct]{#3}}
  \end{picture}}}

\newcommand{\AmidaISmall}[3]
  {\ThreeParallelArcsUpUpUpSmall{}{#1}{#2}{#3}}

\newcommand{\AmidaA}[3]
  {\ThreeParallelArcsUpUpUpSmall{
  \thinlines
  \qbezier[10]( 0.00,10.00)( 2.50,10.00)( 5.00,10.00)
  \thicklines
  }{#1}{#2}{#3}}

\newcommand{\AmidaB}[3]
  {\ThreeParallelArcsUpUpUpSmall{
  \thinlines
  \qbezier[10]( 5.00,10.00)( 7.50,10.00)(10.00,10.00)
  \thicklines
  }{#1}{#2}{#3}}

\newcommand{\AmidaC}[3]
  {\ThreeParallelArcsUpUpUpSmall{
  \thinlines
  \qbezier[20]( 0.00,10.00)( 5.00,10.00)(10.00,10.00)
  \thicklines
  }{#1}{#2}{#3}}

\newcommand{\AmidaI}[3]
  {\ThreeParallelArcsUpUpUp{}{#1}{#2}{#3}}

\newcommand{\AmidaAB}[3]
  {\ThreeParallelArcsUpUpUp{
  \thinlines
  \qbezier[10]( 0.00,10.00)( 2.50,10.00)( 5.00,10.00)
  \qbezier[10]( 5.00, 5.00)( 7.50, 5.00)(10.00, 5.00)
  \thicklines
  }{#1}{#2}{#3}}

\newcommand{\AmidaBA}[3]
  {\ThreeParallelArcsUpUpUp{
  \thinlines
  \qbezier[10]( 5.00,10.00)( 7.50,10.00)(10.00,10.00)
  \qbezier[10]( 0.00, 5.00)( 2.50, 5.00)( 5.00, 5.00)
  \thicklines
  }{#1}{#2}{#3}}

\newcommand{\AmidaBC}[3]
  {\ThreeParallelArcsUpUpUp{
  \thinlines
  \qbezier[10]( 5.00,10.00)( 7.50,10.00)(10.00,10.00)
  \qbezier[20]( 0.00, 5.00)( 5.00, 5.00)(10.00, 5.00)
  \thicklines
  }{#1}{#2}{#3}}

\newcommand{\AmidaCB}[3]
  {\ThreeParallelArcsUpUpUp{
  \thinlines
  \qbezier[20]( 0.00,10.00)( 5.00,10.00)(10.00,10.00)
  \qbezier[10]( 5.00, 5.00)( 7.50, 5.00)(10.00, 5.00)
  \thicklines
  }{#1}{#2}{#3}}

\newcommand{\AmidaCA}[3]
  {\ThreeParallelArcsUpUpUp{
  \thinlines
  \qbezier[20]( 0.00,10.00)( 5.00,10.00)(10.00,10.00)
  \qbezier[10]( 0.00, 5.00)( 2.50, 5.00)( 5.00, 5.00)
  \thicklines
  }{#1}{#2}{#3}}

\newcommand{\AmidaAC}[3]
  {\ThreeParallelArcsUpUpUp{
  \thinlines
  \qbezier[20]( 0.00, 5.00)( 5.00, 5.00)(10.00, 5.00)
  \qbezier[10]( 0.00,10.00)( 2.50,10.00)( 5.00,10.00)
  \thicklines
  }{#1}{#2}{#3}}

\newcommand{\AmidaABA}[3]
  {\ThreeParallelArcsUpUpUp{
  \thinlines
  \qbezier[10]( 0.00, 4.00)( 2.50, 4.00)( 5.00, 4.00)
  \qbezier[10]( 0.00,11.00)( 2.50,11.00)( 5.00,11.00)
  \qbezier[10]( 5.00, 7.50)( 7.50, 7.50)(10.00, 7.50)
  \thicklines
  }{#1}{#2}{#3}}

\newcommand{\AmidaBAB}[3]
  {\ThreeParallelArcsUpUpUp{
  \thinlines
  \qbezier[10]( 5.00, 4.00)( 7.50, 4.00)(10.00, 4.00)
  \qbezier[10]( 5.00,11.00)( 7.50,11.00)(10.00,11.00)
  \qbezier[10]( 0.00, 7.50)( 2.50, 7.50)( 5.00, 7.50)
  \thicklines
  }{#1}{#2}{#3}}

\newcommand{\AmidaIWithCurlTop}[2]
  {\OneArcUpAndOneCurlyArcTop{}{#1}{#2}}

\newcommand{\AmidaABWithCurlTop}[2]
  {\OneArcUpAndOneCurlyArcTop{
  \thinlines
  \qbezier[10]( 0.00,10.00)( 2.50,10.00)( 5.00,10.00)
  \qbezier[10]( 5.00, 5.00)( 7.50, 5.00)(10.00, 5.00)
  \thicklines
  }{#1}{#2}}

\newcommand{\AmidaBAWithCurlTop}[2]
  {\OneArcUpAndOneCurlyArcTop{
  \thinlines
  \qbezier[10]( 5.00,10.00)( 7.50,10.00)(10.00,10.00)
  \qbezier[10]( 0.00, 5.00)( 2.50, 5.00)( 5.00, 5.00)
  \thicklines
  }{#1}{#2}}

\newcommand{\AmidaBCWithCurlTop}[2]
  {\OneArcUpAndOneCurlyArcTop{
  \thinlines
  \qbezier[10]( 5.00,10.00)( 7.50,10.00)(10.00,10.00)
  \qbezier[20]( 0.00, 5.00)( 5.00, 5.00)(10.00, 5.00)
  \thicklines
  }{#1}{#2}}

\newcommand{\AmidaCAWithCurlTop}[2]
  {\OneArcUpAndOneCurlyArcTop{
  \thinlines
  \qbezier[20]( 0.00,10.00)( 5.00,10.00)(10.00,10.00)
  \qbezier[10]( 0.00, 5.00)( 2.50, 5.00)( 5.00, 5.00)
  \thicklines
  }{#1}{#2}}

\newcommand{\AmidaIWithCurlBottom}[2]
  {\OneArcUpAndOneCurlyArcBottom{}{#1}{#2}}

\newcommand{\AmidaABWithCurlBottom}[2]
  {\OneArcUpAndOneCurlyArcBottom{
  \thinlines
  \qbezier[10]( 0.00,10.00)( 2.50,10.00)( 5.00,10.00)
  \qbezier[10]( 5.00, 5.00)( 7.50, 5.00)(10.00, 5.00)
  \thicklines
  }{#1}{#2}}

\newcommand{\AmidaBAWithCurlBottom}[2]
  {\OneArcUpAndOneCurlyArcBottom{
  \thinlines
  \qbezier[10]( 5.00,10.00)( 7.50,10.00)(10.00,10.00)
  \qbezier[10]( 0.00, 5.00)( 2.50, 5.00)( 5.00, 5.00)
  \thicklines
  }{#1}{#2}}

\newcommand{\AmidaBCWithCurlBottom}[2]
  {\OneArcUpAndOneCurlyArcBottom{
  \thinlines
  \qbezier[10]( 5.00,10.00)( 7.50,10.00)(10.00,10.00)
  \qbezier[20]( 0.00, 5.00)( 5.00, 5.00)(10.00, 5.00)
  \thicklines
  }{#1}{#2}}

\newcommand{\AmidaCAWithCurlBottom}[2]
  {\OneArcUpAndOneCurlyArcBottom{
  \thinlines
  \qbezier[20]( 0.00,10.00)( 5.00,10.00)(10.00,10.00)
  \qbezier[10]( 0.00, 5.00)( 2.50, 5.00)( 5.00, 5.00)
  \thicklines
  }{#1}{#2}}

\newcommand{\TwoParallelArcsUpUpWithTwoParallelChords}[2]
  {\TwoParallelArcsUpUp{
  \thinlines
  \qbezier[10]( 0.00, 5.00)( 2.50, 5.00)( 5.00, 5.00)
  \qbezier[10]( 0.00,10.00)( 2.50,10.00)( 5.00,10.00)
  \thicklines
  }{#1}{#2}}

\newcommand{\TwoParallelArcsUpUpWithTwoCrossedChords}[2]
  {\TwoParallelArcsUpUp{
  \thinlines
  \qbezier[14]( 0.00, 5.00)( 2.50, 7.50)( 5.00,10.00)
  \qbezier[14]( 0.00,10.00)( 2.50, 7.50)( 5.00, 5.00)
  \thicklines
  }{#1}{#2}}

\newcommand{\TwoParallelArcsUpUpWithOneChord}[2]
  {\TwoParallelArcsUpUp{
  \thinlines
  \qbezier[10]( 0.00, 7.50)( 2.50, 7.50)( 5.00, 7.50)
  \thicklines
  }{#1}{#2}}

\newcommand{\TwoCrossedArcsUpUp}[2]
  {\raisebox{-6.75mm}
  {\begin{picture}( 5.00,15.00)
  \put( 0.00, 0.00){\vector( 1, 3){ 5.00}}
  \put( 5.00, 0.00){\vector(-1, 3){ 5.00}}
  \put( 0.00,-0.50){\makebox( 0.00, 0.00)[ct]{#1}}
  \put( 5.00,-0.50){\makebox( 0.00, 0.00)[ct]{#2}}
  \end{picture}}}

\newcommand{\DoublePointUpUp}[2]
  {\raisebox{-6.75mm}
  {\begin{picture}( 5.00,15.00)
  \put( 0.00, 0.00){\vector( 1, 3){ 5.00}}
  \put( 5.00, 0.00){\vector(-1, 3){ 5.00}}
  \put( 2.50, 7.50){\circle*{2.00}}
  \put( 0.00,-0.50){\makebox( 0.00, 0.00)[ct]{#1}}
  \put( 5.00,-0.50){\makebox( 0.00, 0.00)[ct]{#2}}
  \end{picture}}}

\newcommand{\TwoCrossedArcsUpUpWithOneChord}[2]
  {\raisebox{-6.75mm}
  {\begin{picture}( 5.00,15.00)
  \put( 0.00, 0.00){\vector( 1, 3){ 5.00}}
  \put( 5.00, 0.00){\vector(-1, 3){ 5.00}}
  \thinlines
  \qbezier[6]( 1.00, 3.00)( 2.50, 3.00)(4.00, 3.00)
  \thicklines
  \put( 0.00,-0.50){\makebox( 0.00, 0.00)[ct]{#1}}
  \put( 5.00,-0.50){\makebox( 0.00, 0.00)[ct]{#2}}
  \end{picture}}}

\newcommand{\OneArcUpAndOneCircleWithOneChord}[2]
  {\raisebox{-6.75mm}
  {\begin{picture}(15.00,15.00)
  \put( 0.00, 0.00){\vector( 0, 1){16.00}}
  \put(10.00, 7.50){\circle{10.00}}
  \put(15.00, 8.50){\vector( 0, 1){ 0.00}}
  \thinlines
  \qbezier[10]( 0.00, 7.50)( 2.50, 7.50)( 5.00, 7.50)
  \thicklines
  \put( 0.00,-0.50){\makebox( 0.00, 0.00)[ct]{#1}}
  \put(10.00, 2.00){\makebox( 0.00, 0.00)[ct]{#2}}
  \end{picture}}}

\newcommand{\TwoArcsUpUpAndOneCircleWithOneChordToCircle}[3]
  {\raisebox{-6.75mm}
  {\begin{picture}(15.00,15.00)
  \put( 0.00, 0.00){\vector( 0, 1){16.00}}
  \put(15.00, 0.00){\vector( 0, 1){16.00}} 
  \put( 7.50, 7.50){\circle{ 5.00}}
  \put(10.00, 8.50){\vector( 0, 1){ 0.00}}
  \thinlines
  \qbezier[10]( 0.00, 7.50)( 2.50, 7.50)( 5.00, 7.50)
  \thicklines
  \put( 0.00,-0.50){\makebox( 0.00, 0.00)[ct]{#1}}
  \put( 7.50, 4.50){\makebox( 0.00, 0.00)[ct]{#2}}
  \put(15.00,-0.50){\makebox( 0.00, 0.00)[ct]{#3}}
  \end{picture}}}

\newcommand{\TwoArcsUpUpAndOneCircleWithOneChordToArcs}[3]
  {\raisebox{-6.75mm}
  {\begin{picture}(15.00,15.00)
  \put( 0.00, 0.00){\vector( 0, 1){16.00}}
  \put( 5.00, 0.00){\vector( 0, 1){16.00}} 
  \put(12.50, 7.50){\circle{ 5.00}}
  \put(15.00, 8.50){\vector( 0, 1){ 0.00}}
  \thinlines
  \qbezier[10]( 0.00, 7.50)( 2.50, 7.50)( 5.00, 7.50)
  \thicklines
  \put( 0.00,-0.50){\makebox( 0.00, 0.00)[ct]{#1}}
  \put( 5.00,-0.50){\makebox( 0.00, 0.00)[ct]{#2}}
  \put(15.00, 4.50){\makebox( 0.00, 0.00)[ct]{#3}}
  \end{picture}}}

\newcommand{\OneArcUpAndOneCircleWithoutChord}[2]
  {\raisebox{-6.75mm}
  {\begin{picture}(15.00,15.00)
  \put( 0.00, 0.00){\vector( 0, 1){16.00}}
  \put(10.00, 7.50){\circle{10.00}}
  \put(15.00, 8.50){\vector( 0, 1){ 0.00}}
  \put( 0.00,-0.50){\makebox( 0.00, 0.00)[ct]{#1}}
  \put(10.00, 2.00){\makebox( 0.00, 0.00)[ct]{#2}}
  \end{picture}}}

\newcommand{\TwoCirclesWithOneChord}[2]
  {\raisebox{-4.50mm}
  {\begin{picture}(25.00,10.00)
  \put( 5.00, 5.00){\circle{10.00}}
  \put( 6.00,10.00){\vector( 1, 0){ 0.00}}
  \put(20.00, 5.00){\circle{10.00}}
  \put(21.00,10.00){\vector( 1, 0){ 0.00}}
  \thinlines
  \qbezier[10](10.00, 5.00)(12.50, 5.00)(15.00, 5.00)
  \thicklines
  \put( 5.00,-0.50){\makebox( 0.00, 0.00)[ct]{#1}}
  \put(20.00,-0.50){\makebox( 0.00, 0.00)[ct]{#2}}
  \end{picture}}}

\newcommand{\OneArcUp}[1]
  {\raisebox{-6.75mm}
  {\begin{picture}( 5.00,15.00)
  \put( 2.50, 0.00){\vector( 0, 1){16.00}}
  \put( 2.50,-0.50){\makebox( 0.00, 0.00)[ct]{#1}}
  \end{picture}}}

\newcommand{\OneCircle}[1]
  {\raisebox{-4.50mm}
  {\begin{picture}(10.00,10.00)
  \put( 5.00, 5.00){\circle{10.00}}
  \put(10.00, 6.00){\vector( 0, 1){ 0.00}}
  \put( 5.00,-0.50){\makebox( 0.00, 0.00)[ct]{#1}}
  \end{picture}}}

\newcommand{\PositiveCrossUpUp}[2]
  {\raisebox{-6.75mm}
  {\begin{picture}( 5.00,15.00)
  \qbezier( 0.00, 5.00)( 0.00, 6.00)( 2.50, 7.50)
  \qbezier( 5.00,10.00)( 5.00, 9.00)( 2.50, 7.50)
  \qbezier( 5.00, 5.00)( 5.00, 6.00)( 3.50, 7.00)
  \qbezier( 0.00,10.00)( 0.00, 9.00)( 1.50, 8.00)
  \put( 0.00,10.00){\vector( 0, 1){ 6.00}}
  \put( 5.00,10.00){\vector( 0, 1){ 6.00}}
  \put( 0.00, 0.00){\line( 0, 1){ 5.00}}
  \put( 5.00, 0.00){\line( 0, 1){ 5.00}}
  \put( 0.00,-0.50){\makebox( 0.00, 0.00)[ct]{#1}}
  \put( 5.00,-0.50){\makebox( 0.00, 0.00)[ct]{#2}}
  \end{picture}}}

\newcommand{\NegativeCrossUpUp}[2]
  {\raisebox{-6.75mm}
  {\begin{picture}( 5.00,15.00)
  \qbezier( 0.00, 5.00)( 0.00, 6.00)( 1.50, 7.00)
  \qbezier( 5.00,10.00)( 5.00, 9.00)( 3.50, 8.00)
  \qbezier( 5.00, 5.00)( 5.00, 6.00)( 2.50, 7.50)
  \qbezier( 0.00,10.00)( 0.00, 9.00)( 2.50, 7.50)
  \put( 0.00,10.00){\vector( 0, 1){ 6.00}}
  \put( 5.00,10.00){\vector( 0, 1){ 6.00}}
  \put( 0.00, 0.00){\line( 0, 1){ 5.00}}
  \put( 5.00, 0.00){\line( 0, 1){ 5.00}}
  \put( 0.00,-0.50){\makebox( 0.00, 0.00)[ct]{#1}}
  \put( 5.00,-0.50){\makebox( 0.00, 0.00)[ct]{#2}}
  \end{picture}}}

\newcommand{\PositiveFullTwistUpUp}[2]
  {\raisebox{-6.75mm}
  {\begin{picture}( 5.00,15.00)
  \qbezier( 0.00, 2.50)( 0.00, 3.50)( 2.50, 5.00)
  \qbezier( 5.00, 7.50)( 5.00, 6.50)( 2.50, 5.00)
  \qbezier( 5.00, 2.50)( 5.00, 3.50)( 3.50, 4.50)
  \qbezier( 0.00, 7.50)( 0.00, 6.50)( 1.50, 5.50)
  \qbezier( 0.00, 7.50)( 0.00, 8.50)( 2.50,10.00)
  \qbezier( 5.00,12.50)( 5.00,11.50)( 2.50,10.00)
  \qbezier( 5.00, 7.50)( 5.00, 8.50)( 3.50, 9.50)
  \qbezier( 0.00,12.50)( 0.00,11.50)( 1.50,10.50)
  \put( 0.00,12.50){\vector( 0, 1){ 3.50}}
  \put( 5.00,12.50){\vector( 0, 1){ 3.50}}
  \put( 0.00, 0.00){\line( 0, 1){ 2.50}}
  \put( 5.00, 0.00){\line( 0, 1){ 2.50}}
  \put( 0.00,-0.50){\makebox( 0.00, 0.00)[ct]{#1}}
  \put( 5.00,-0.50){\makebox( 0.00, 0.00)[ct]{#2}}
  \end{picture}}}

\newcommand{\TwoBraidCP}[2]
  {\raisebox{-6.75mm}
  {\begin{picture}( 5.00,15.00)
  \qbezier( 0.00, 2.50)( 0.00, 3.50)( 2.50, 5.00)
  \qbezier( 5.00, 7.50)( 5.00, 6.50)( 2.50, 5.00)
  \qbezier( 5.00, 2.50)( 5.00, 3.50)( 3.50, 4.50)
  \qbezier( 0.00, 7.50)( 0.00, 6.50)( 1.50, 5.50)
  \qbezier( 0.00, 7.50)( 0.00, 8.50)( 2.50,10.00)
  \qbezier( 5.00,12.50)( 5.00,11.50)( 2.50,10.00)
  \qbezier( 5.00, 7.50)( 5.00, 8.50)( 2.50,10.00)
  \qbezier( 0.00,12.50)( 0.00,11.50)( 2.50,10.00)
  \put( 0.00,12.50){\vector( 0, 1){ 3.50}}
  \put( 5.00,12.50){\vector( 0, 1){ 3.50}}
  \put( 0.00, 0.00){\line( 0, 1){ 2.50}}
  \put( 5.00, 0.00){\line( 0, 1){ 2.50}}
  \put( 2.50,10.00){\circle*{ 2.00}}
  \put( 0.00,-0.50){\makebox( 0.00, 0.00)[ct]{#1}}
  \put( 5.00,-0.50){\makebox( 0.00, 0.00)[ct]{#2}}
  \end{picture}}}

\newcommand{\NegativeFullTwistUpUp}[2]
  {\raisebox{-6.75mm}
  {\begin{picture}( 5.00,15.00)
  \qbezier( 0.00, 2.50)( 0.00, 3.50)( 1.50, 4.50)
  \qbezier( 5.00, 7.50)( 5.00, 6.50)( 3.50, 5.50)
  \qbezier( 5.00, 2.50)( 5.00, 3.50)( 2.50, 5.00)
  \qbezier( 0.00, 7.50)( 0.00, 6.50)( 2.50, 5.00)
  \qbezier( 0.00, 7.50)( 0.00, 8.50)( 1.50, 9.50)
  \qbezier( 5.00,12.50)( 5.00,11.50)( 3.50,10.50)
  \qbezier( 5.00, 7.50)( 5.00, 8.50)( 2.50,10.00)
  \qbezier( 0.00,12.50)( 0.00,11.50)( 2.50,10.00)
  \put( 0.00,12.50){\vector( 0, 1){ 3.50}}
  \put( 5.00,12.50){\vector( 0, 1){ 3.50}}
  \put( 0.00, 0.00){\line( 0, 1){ 2.50}}
  \put( 5.00, 0.00){\line( 0, 1){ 2.50}}
  \put( 0.00,-0.50){\makebox( 0.00, 0.00)[ct]{#1}}
  \put( 5.00,-0.50){\makebox( 0.00, 0.00)[ct]{#2}}
  \end{picture}}}

\newcommand{\ThreeBraidPOnePTwoPOne}[3]
  {\raisebox{-6.75mm}
  {\begin{picture}(10.00,15.00)
  \qbezier( 0.00, 0.00)( 0.00, 1.00)( 2.50, 2.50)
  \qbezier( 5.00, 5.00)( 5.00, 4.00)( 2.50, 2.50)
  \qbezier( 5.00, 0.00)( 5.00, 1.00)( 3.50, 2.00)
  \qbezier( 0.00, 5.00)( 0.00, 4.00)( 1.50, 3.00)
  \qbezier( 5.00, 5.00)( 5.00, 6.00)( 7.50, 7.50)
  \qbezier(10.00,10.00)(10.00, 9.00)( 7.50, 7.50)
  \qbezier(10.00, 5.00)(10.00, 6.00)( 8.50, 7.00)
  \qbezier( 5.00,10.00)( 5.00, 9.00)( 6.50, 8.00)
  \qbezier( 0.00,10.00)( 0.00,11.00)( 2.50,12.50)
  \qbezier( 5.00,15.00)( 5.00,14.00)( 2.50,12.50)
  \qbezier( 5.00,10.00)( 5.00,11.00)( 3.50,12.00)
  \qbezier( 0.00,15.00)( 0.00,14.00)( 1.50,13.00)
  \put( 0.00, 5.00){\line( 0, 1){ 5.00}}
  \put(10.00, 0.00){\line( 0, 1){ 5.00}}
  \put( 0.00,15.00){\vector( 0, 1){ 2.00}}
  \put( 5.00,15.00){\vector( 0, 1){ 2.00}}
  \put(10.00,10.00){\vector( 0, 1){ 7.00}}
  \put( 0.00,-0.50){\makebox( 0.00, 0.00)[ct]{#1}}
  \put( 5.00,-0.50){\makebox( 0.00, 0.00)[ct]{#2}}
  \put(10.00,-0.50){\makebox( 0.00, 0.00)[ct]{#3}}
  \end{picture}}}

\newcommand{\ThreeBraidPOnePTwoNOne}[3]
  {\raisebox{-6.75mm}
  {\begin{picture}(10.00,15.00)
  \qbezier( 0.00, 0.00)( 0.00, 1.00)( 2.50, 2.50)
  \qbezier( 5.00, 5.00)( 5.00, 4.00)( 2.50, 2.50)
  \qbezier( 5.00, 0.00)( 5.00, 1.00)( 3.50, 2.00)
  \qbezier( 0.00, 5.00)( 0.00, 4.00)( 1.50, 3.00)
  \qbezier( 5.00, 5.00)( 5.00, 6.00)( 7.50, 7.50)
  \qbezier(10.00,10.00)(10.00, 9.00)( 7.50, 7.50)
  \qbezier(10.00, 5.00)(10.00, 6.00)( 8.50, 7.00)
  \qbezier( 5.00,10.00)( 5.00, 9.00)( 6.50, 8.00)
  \qbezier( 0.00,10.00)( 0.00,11.00)( 1.50,12.00)
  \qbezier( 5.00,15.00)( 5.00,14.00)( 3.50,13.00)
  \qbezier( 5.00,10.00)( 5.00,11.00)( 2.50,12.50)
  \qbezier( 0.00,15.00)( 0.00,14.00)( 2.50,12.50)
  \put( 0.00, 5.00){\line( 0, 1){ 5.00}}
  \put(10.00, 0.00){\line( 0, 1){ 5.00}}
  \put( 0.00,15.00){\vector( 0, 1){ 2.00}}
  \put( 5.00,15.00){\vector( 0, 1){ 2.00}}
  \put(10.00,10.00){\vector( 0, 1){ 7.00}}
  \put( 0.00,-0.50){\makebox( 0.00, 0.00)[ct]{#1}}
  \put( 5.00,-0.50){\makebox( 0.00, 0.00)[ct]{#2}}
  \put(10.00,-0.50){\makebox( 0.00, 0.00)[ct]{#3}}
  \end{picture}}}

\newcommand{\ThreeBraidPOneNTwoPOne}[3]
  {\raisebox{-6.75mm}
  {\begin{picture}(10.00,15.00)
  \qbezier( 0.00, 0.00)( 0.00, 1.00)( 2.50, 2.50)
  \qbezier( 5.00, 5.00)( 5.00, 4.00)( 2.50, 2.50)
  \qbezier( 5.00, 0.00)( 5.00, 1.00)( 3.50, 2.00)
  \qbezier( 0.00, 5.00)( 0.00, 4.00)( 1.50, 3.00)
  \qbezier( 5.00, 5.00)( 5.00, 6.00)( 6.50, 7.00)
  \qbezier(10.00,10.00)(10.00, 9.00)( 8.50, 8.00)
  \qbezier(10.00, 5.00)(10.00, 6.00)( 7.50, 7.50)
  \qbezier( 5.00,10.00)( 5.00, 9.00)( 7.50, 7.50)
  \qbezier( 0.00,10.00)( 0.00,11.00)( 2.50,12.50)
  \qbezier( 5.00,15.00)( 5.00,14.00)( 2.50,12.50)
  \qbezier( 5.00,10.00)( 5.00,11.00)( 3.50,12.00)
  \qbezier( 0.00,15.00)( 0.00,14.00)( 1.50,13.00)
  \put( 0.00, 5.00){\line( 0, 1){ 5.00}}
  \put(10.00, 0.00){\line( 0, 1){ 5.00}}
  \put( 0.00,15.00){\vector( 0, 1){ 2.00}}
  \put( 5.00,15.00){\vector( 0, 1){ 2.00}}
  \put(10.00,10.00){\vector( 0, 1){ 7.00}}
  \put( 0.00,-0.50){\makebox( 0.00, 0.00)[ct]{#1}}
  \put( 5.00,-0.50){\makebox( 0.00, 0.00)[ct]{#2}}
  \put(10.00,-0.50){\makebox( 0.00, 0.00)[ct]{#3}}
  \end{picture}}}

\newcommand{\ThreeBraidNOnePTwoPOne}[3]
  {\raisebox{-6.75mm}
  {\begin{picture}(10.00,15.00)
  \qbezier( 0.00, 0.00)( 0.00, 1.00)( 1.50, 2.00)
  \qbezier( 5.00, 5.00)( 5.00, 4.00)( 3.50, 3.00)
  \qbezier( 5.00, 0.00)( 5.00, 1.00)( 2.50, 2.50)
  \qbezier( 0.00, 5.00)( 0.00, 4.00)( 2.50, 2.50)
  \qbezier( 5.00, 5.00)( 5.00, 6.00)( 7.50, 7.50)
  \qbezier(10.00,10.00)(10.00, 9.00)( 7.50, 7.50)
  \qbezier(10.00, 5.00)(10.00, 6.00)( 8.50, 7.00)
  \qbezier( 5.00,10.00)( 5.00, 9.00)( 6.50, 8.00)
  \qbezier( 0.00,10.00)( 0.00,11.00)( 2.50,12.50)
  \qbezier( 5.00,15.00)( 5.00,14.00)( 2.50,12.50)
  \qbezier( 5.00,10.00)( 5.00,11.00)( 3.50,12.00)
  \qbezier( 0.00,15.00)( 0.00,14.00)( 1.50,13.00)
  \put( 0.00, 5.00){\line( 0, 1){ 5.00}}
  \put(10.00, 0.00){\line( 0, 1){ 5.00}}
  \put( 0.00,15.00){\vector( 0, 1){ 2.00}}
  \put( 5.00,15.00){\vector( 0, 1){ 2.00}}
  \put(10.00,10.00){\vector( 0, 1){ 7.00}}
  \put( 0.00,-0.50){\makebox( 0.00, 0.00)[ct]{#1}}
  \put( 5.00,-0.50){\makebox( 0.00, 0.00)[ct]{#2}}
  \put(10.00,-0.50){\makebox( 0.00, 0.00)[ct]{#3}}
  \end{picture}}}

\newcommand{\ThreeBraidPOneNTwoNOne}[3]
  {\raisebox{-6.75mm}
  {\begin{picture}(10.00,15.00)
  \qbezier( 0.00, 0.00)( 0.00, 1.00)( 2.50, 2.50)
  \qbezier( 5.00, 5.00)( 5.00, 4.00)( 2.50, 2.50)
  \qbezier( 5.00, 0.00)( 5.00, 1.00)( 3.50, 2.00)
  \qbezier( 0.00, 5.00)( 0.00, 4.00)( 1.50, 3.00)
  \qbezier( 5.00, 5.00)( 5.00, 6.00)( 6.50, 7.00)
  \qbezier(10.00,10.00)(10.00, 9.00)( 8.50, 8.00)
  \qbezier(10.00, 5.00)(10.00, 6.00)( 7.50, 7.50)
  \qbezier( 5.00,10.00)( 5.00, 9.00)( 7.50, 7.50)
  \qbezier( 0.00,10.00)( 0.00,11.00)( 1.50,12.00)
  \qbezier( 5.00,15.00)( 5.00,14.00)( 3.50,13.00)
  \qbezier( 5.00,10.00)( 5.00,11.00)( 2.50,12.50)
  \qbezier( 0.00,15.00)( 0.00,14.00)( 2.50,12.50)
  \put( 0.00, 5.00){\line( 0, 1){ 5.00}}
  \put(10.00, 0.00){\line( 0, 1){ 5.00}}
  \put( 0.00,15.00){\vector( 0, 1){ 2.00}}
  \put( 5.00,15.00){\vector( 0, 1){ 2.00}}
  \put(10.00,10.00){\vector( 0, 1){ 7.00}}
  \put( 0.00,-0.50){\makebox( 0.00, 0.00)[ct]{#1}}
  \put( 5.00,-0.50){\makebox( 0.00, 0.00)[ct]{#2}}
  \put(10.00,-0.50){\makebox( 0.00, 0.00)[ct]{#3}}
  \end{picture}}}

\newcommand{\ThreeBraidNOnePTwoNOne}[3]
  {\raisebox{-6.75mm}
  {\begin{picture}(10.00,15.00)
  \qbezier( 0.00, 0.00)( 0.00, 1.00)( 1.50, 2.00)
  \qbezier( 5.00, 5.00)( 5.00, 4.00)( 3.50, 3.00)
  \qbezier( 5.00, 0.00)( 5.00, 1.00)( 2.50, 2.50)
  \qbezier( 0.00, 5.00)( 0.00, 4.00)( 2.50, 2.50)
  \qbezier( 5.00, 5.00)( 5.00, 6.00)( 7.50, 7.50)
  \qbezier(10.00,10.00)(10.00, 9.00)( 7.50, 7.50)
  \qbezier(10.00, 5.00)(10.00, 6.00)( 8.50, 7.00)
  \qbezier( 5.00,10.00)( 5.00, 9.00)( 6.50, 8.00)
  \qbezier( 0.00,10.00)( 0.00,11.00)( 1.50,12.00)
  \qbezier( 5.00,15.00)( 5.00,14.00)( 3.50,13.00)
  \qbezier( 5.00,10.00)( 5.00,11.00)( 2.50,12.50)
  \qbezier( 0.00,15.00)( 0.00,14.00)( 2.50,12.50)
  \put( 0.00, 5.00){\line( 0, 1){ 5.00}}
  \put(10.00, 0.00){\line( 0, 1){ 5.00}}
  \put( 0.00,15.00){\vector( 0, 1){ 2.00}}
  \put( 5.00,15.00){\vector( 0, 1){ 2.00}}
  \put(10.00,10.00){\vector( 0, 1){ 7.00}}
  \put( 0.00,-0.50){\makebox( 0.00, 0.00)[ct]{#1}}
  \put( 5.00,-0.50){\makebox( 0.00, 0.00)[ct]{#2}}
  \put(10.00,-0.50){\makebox( 0.00, 0.00)[ct]{#3}}
  \end{picture}}}

\newcommand{\ThreeBraidNOneNTwoPOne}[3]
  {\raisebox{-6.75mm}
  {\begin{picture}(10.00,15.00)
  \qbezier( 0.00, 0.00)( 0.00, 1.00)( 1.50, 2.00)
  \qbezier( 5.00, 5.00)( 5.00, 4.00)( 3.50, 3.00)
  \qbezier( 5.00, 0.00)( 5.00, 1.00)( 2.50, 2.50)
  \qbezier( 0.00, 5.00)( 0.00, 4.00)( 2.50, 2.50)
  \qbezier( 5.00, 5.00)( 5.00, 6.00)( 6.50, 7.00)
  \qbezier(10.00,10.00)(10.00, 9.00)( 8.50, 8.00)
  \qbezier(10.00, 5.00)(10.00, 6.00)( 7.50, 7.50)
  \qbezier( 5.00,10.00)( 5.00, 9.00)( 7.50, 7.50)
  \qbezier( 0.00,10.00)( 0.00,11.00)( 2.50,12.50)
  \qbezier( 5.00,15.00)( 5.00,14.00)( 2.50,12.50)
  \qbezier( 5.00,10.00)( 5.00,11.00)( 3.50,12.00)
  \qbezier( 0.00,15.00)( 0.00,14.00)( 1.50,13.00)
  \put( 0.00, 5.00){\line( 0, 1){ 5.00}}
  \put(10.00, 0.00){\line( 0, 1){ 5.00}}
  \put( 0.00,15.00){\vector( 0, 1){ 2.00}}
  \put( 5.00,15.00){\vector( 0, 1){ 2.00}}
  \put(10.00,10.00){\vector( 0, 1){ 7.00}}
  \put( 0.00,-0.50){\makebox( 0.00, 0.00)[ct]{#1}}
  \put( 5.00,-0.50){\makebox( 0.00, 0.00)[ct]{#2}}
  \put(10.00,-0.50){\makebox( 0.00, 0.00)[ct]{#3}}
  \end{picture}}}

\newcommand{\ThreeBraidNOneNTwoNOne}[3]
  {\raisebox{-6.75mm}
  {\begin{picture}(10.00,15.00)
  \qbezier( 0.00, 0.00)( 0.00, 1.00)( 1.50, 2.00)
  \qbezier( 5.00, 5.00)( 5.00, 4.00)( 3.50, 3.00)
  \qbezier( 5.00, 0.00)( 5.00, 1.00)( 2.50, 2.50)
  \qbezier( 0.00, 5.00)( 0.00, 4.00)( 2.50, 2.50)
  \qbezier( 5.00, 5.00)( 5.00, 6.00)( 6.50, 7.00)
  \qbezier(10.00,10.00)(10.00, 9.00)( 8.50, 8.00)
  \qbezier(10.00, 5.00)(10.00, 6.00)( 7.50, 7.50)
  \qbezier( 5.00,10.00)( 5.00, 9.00)( 7.50, 7.50)
  \qbezier( 0.00,10.00)( 0.00,11.00)( 1.50,12.00)
  \qbezier( 5.00,15.00)( 5.00,14.00)( 3.50,13.00)
  \qbezier( 5.00,10.00)( 5.00,11.00)( 2.50,12.50)
  \qbezier( 0.00,15.00)( 0.00,14.00)( 2.50,12.50)
  \put( 0.00, 5.00){\line( 0, 1){ 5.00}}
  \put(10.00, 0.00){\line( 0, 1){ 5.00}}
  \put( 0.00,15.00){\vector( 0, 1){ 2.00}}
  \put( 5.00,15.00){\vector( 0, 1){ 2.00}}
  \put(10.00,10.00){\vector( 0, 1){ 7.00}}
  \put( 0.00,-0.50){\makebox( 0.00, 0.00)[ct]{#1}}
  \put( 5.00,-0.50){\makebox( 0.00, 0.00)[ct]{#2}}
  \put(10.00,-0.50){\makebox( 0.00, 0.00)[ct]{#3}}
  \end{picture}}}

\newcommand{\ThreeBraidPTwoPOnePTwo}[3]
  {\raisebox{-6.75mm}
  {\begin{picture}(10.00,15.00)
  \qbezier( 5.00, 0.00)( 5.00, 1.00)( 7.50, 2.50)
  \qbezier(10.00, 5.00)(10.00, 4.00)( 7.50, 2.50)
  \qbezier(10.00, 0.00)(10.00, 1.00)( 8.50, 2.00)
  \qbezier( 5.00, 5.00)( 5.00, 4.00)( 6.50, 3.00)
  \qbezier( 0.00, 5.00)( 0.00, 6.00)( 2.50, 7.50)
  \qbezier( 5.00,10.00)( 5.00, 9.00)( 2.50, 7.50)
  \qbezier( 5.00, 5.00)( 5.00, 6.00)( 3.50, 7.00)
  \qbezier( 0.00,10.00)( 0.00, 9.00)( 1.50, 8.00)
  \qbezier( 5.00,10.00)( 5.00,11.00)( 7.50,12.50)
  \qbezier(10.00,15.00)(10.00,14.00)( 7.50,12.50)
  \qbezier(10.00,10.00)(10.00,11.00)( 8.50,12.00)
  \qbezier( 5.00,15.00)( 5.00,14.00)( 6.50,13.00)
  \put(10.00, 5.00){\line( 0, 1){ 5.00}}
  \put( 0.00, 0.00){\line( 0, 1){ 5.00}}
  \put(10.00,15.00){\vector( 0, 1){ 2.00}}
  \put( 5.00,15.00){\vector( 0, 1){ 2.00}}
  \put( 0.00,10.00){\vector( 0, 1){ 7.00}}
  \put( 0.00,-0.50){\makebox( 0.00, 0.00)[ct]{#1}}
  \put( 5.00,-0.50){\makebox( 0.00, 0.00)[ct]{#2}}
  \put(10.00,-0.50){\makebox( 0.00, 0.00)[ct]{#3}}
  \end{picture}}}

\newcommand{\ThreeBraidPTwoPOneNTwo}[3]
  {\raisebox{-6.75mm}
  {\begin{picture}(10.00,15.00)
  \qbezier( 5.00, 0.00)( 5.00, 1.00)( 7.50, 2.50)
  \qbezier(10.00, 5.00)(10.00, 4.00)( 7.50, 2.50)
  \qbezier(10.00, 0.00)(10.00, 1.00)( 8.50, 2.00)
  \qbezier( 5.00, 5.00)( 5.00, 4.00)( 6.50, 3.00)
  \qbezier( 0.00, 5.00)( 0.00, 6.00)( 2.50, 7.50)
  \qbezier( 5.00,10.00)( 5.00, 9.00)( 2.50, 7.50)
  \qbezier( 5.00, 5.00)( 5.00, 6.00)( 3.50, 7.00)
  \qbezier( 0.00,10.00)( 0.00, 9.00)( 1.50, 8.00)
  \qbezier( 5.00,10.00)( 5.00,11.00)( 6.50,12.00)
  \qbezier(10.00,15.00)(10.00,14.00)( 8.50,13.00)
  \qbezier(10.00,10.00)(10.00,11.00)( 7.50,12.50)
  \qbezier( 5.00,15.00)( 5.00,14.00)( 7.50,12.50)
  \put(10.00, 5.00){\line( 0, 1){ 5.00}}
  \put( 0.00, 0.00){\line( 0, 1){ 5.00}}
  \put(10.00,15.00){\vector( 0, 1){ 2.00}}
  \put( 5.00,15.00){\vector( 0, 1){ 2.00}}
  \put( 0.00,10.00){\vector( 0, 1){ 7.00}}
  \put( 0.00,-0.50){\makebox( 0.00, 0.00)[ct]{#1}}
  \put( 5.00,-0.50){\makebox( 0.00, 0.00)[ct]{#2}}
  \put(10.00,-0.50){\makebox( 0.00, 0.00)[ct]{#3}}
  \end{picture}}}

\newcommand{\ThreeBraidPTwoNOnePTwo}[3]
  {\raisebox{-6.75mm}
  {\begin{picture}(10.00,15.00)
  \qbezier( 5.00, 0.00)( 5.00, 1.00)( 7.50, 2.50)
  \qbezier(10.00, 5.00)(10.00, 4.00)( 7.50, 2.50)
  \qbezier(10.00, 0.00)(10.00, 1.00)( 8.50, 2.00)
  \qbezier( 5.00, 5.00)( 5.00, 4.00)( 6.50, 3.00)
  \qbezier( 0.00, 5.00)( 0.00, 6.00)( 1.50, 7.00)
  \qbezier( 5.00,10.00)( 5.00, 9.00)( 3.50, 8.00)
  \qbezier( 5.00, 5.00)( 5.00, 6.00)( 2.50, 7.50)
  \qbezier( 0.00,10.00)( 0.00, 9.00)( 2.50, 7.50)
  \qbezier( 5.00,10.00)( 5.00,11.00)( 7.50,12.50)
  \qbezier(10.00,15.00)(10.00,14.00)( 7.50,12.50)
  \qbezier(10.00,10.00)(10.00,11.00)( 8.50,12.00)
  \qbezier( 5.00,15.00)( 5.00,14.00)( 6.50,13.00)
  \put(10.00, 5.00){\line( 0, 1){ 5.00}}
  \put( 0.00, 0.00){\line( 0, 1){ 5.00}}
  \put(10.00,15.00){\vector( 0, 1){ 2.00}}
  \put( 5.00,15.00){\vector( 0, 1){ 2.00}}
  \put( 0.00,10.00){\vector( 0, 1){ 7.00}}
  \put( 0.00,-0.50){\makebox( 0.00, 0.00)[ct]{#1}}
  \put( 5.00,-0.50){\makebox( 0.00, 0.00)[ct]{#2}}
  \put(10.00,-0.50){\makebox( 0.00, 0.00)[ct]{#3}}
  \end{picture}}}

\newcommand{\ThreeBraidNTwoPOnePTwo}[3]
  {\raisebox{-6.75mm}
  {\begin{picture}(10.00,15.00)
  \qbezier( 5.00, 0.00)( 5.00, 1.00)( 6.50, 2.00)
  \qbezier(10.00, 5.00)(10.00, 4.00)( 8.50, 3.00)
  \qbezier(10.00, 0.00)(10.00, 1.00)( 7.50, 2.50)
  \qbezier( 5.00, 5.00)( 5.00, 4.00)( 7.50, 2.50)
  \qbezier( 0.00, 5.00)( 0.00, 6.00)( 2.50, 7.50)
  \qbezier( 5.00,10.00)( 5.00, 9.00)( 2.50, 7.50)
  \qbezier( 5.00, 5.00)( 5.00, 6.00)( 3.50, 7.00)
  \qbezier( 0.00,10.00)( 0.00, 9.00)( 1.50, 8.00)
  \qbezier( 5.00,10.00)( 5.00,11.00)( 7.50,12.50)
  \qbezier(10.00,15.00)(10.00,14.00)( 7.50,12.50)
  \qbezier(10.00,10.00)(10.00,11.00)( 8.50,12.00)
  \qbezier( 5.00,15.00)( 5.00,14.00)( 6.50,13.00)
  \put(10.00, 5.00){\line( 0, 1){ 5.00}}
  \put( 0.00, 0.00){\line( 0, 1){ 5.00}}
  \put(10.00,15.00){\vector( 0, 1){ 2.00}}
  \put( 5.00,15.00){\vector( 0, 1){ 2.00}}
  \put( 0.00,10.00){\vector( 0, 1){ 7.00}}
  \put( 0.00,-0.50){\makebox( 0.00, 0.00)[ct]{#1}}
  \put( 5.00,-0.50){\makebox( 0.00, 0.00)[ct]{#2}}
  \put(10.00,-0.50){\makebox( 0.00, 0.00)[ct]{#3}}
  \end{picture}}}

\newcommand{\ThreeBraidPTwoNOneNTwo}[3]
  {\raisebox{-6.75mm}
  {\begin{picture}(10.00,15.00)
  \qbezier( 5.00, 0.00)( 5.00, 1.00)( 7.50, 2.50)
  \qbezier(10.00, 5.00)(10.00, 4.00)( 7.50, 2.50)
  \qbezier(10.00, 0.00)(10.00, 1.00)( 8.50, 2.00)
  \qbezier( 5.00, 5.00)( 5.00, 4.00)( 6.50, 3.00)
  \qbezier( 0.00, 5.00)( 0.00, 6.00)( 1.50, 7.00)
  \qbezier( 5.00,10.00)( 5.00, 9.00)( 3.50, 8.00)
  \qbezier( 5.00, 5.00)( 5.00, 6.00)( 2.50, 7.50)
  \qbezier( 0.00,10.00)( 0.00, 9.00)( 2.50, 7.50)
  \qbezier( 5.00,10.00)( 5.00,11.00)( 6.50,12.00)
  \qbezier(10.00,15.00)(10.00,14.00)( 8.50,13.00)
  \qbezier(10.00,10.00)(10.00,11.00)( 7.50,12.50)
  \qbezier( 5.00,15.00)( 5.00,14.00)( 7.50,12.50)
  \put(10.00, 5.00){\line( 0, 1){ 5.00}}
  \put( 0.00, 0.00){\line( 0, 1){ 5.00}}
  \put(10.00,15.00){\vector( 0, 1){ 2.00}}
  \put( 5.00,15.00){\vector( 0, 1){ 2.00}}
  \put( 0.00,10.00){\vector( 0, 1){ 7.00}}
  \put( 0.00,-0.50){\makebox( 0.00, 0.00)[ct]{#1}}
  \put( 5.00,-0.50){\makebox( 0.00, 0.00)[ct]{#2}}
  \put(10.00,-0.50){\makebox( 0.00, 0.00)[ct]{#3}}
  \end{picture}}}

\newcommand{\ThreeBraidNTwoPOneNTwo}[3]
  {\raisebox{-6.75mm}
  {\begin{picture}(10.00,15.00)
  \qbezier( 5.00, 0.00)( 5.00, 1.00)( 6.50, 2.00)
  \qbezier(10.00, 5.00)(10.00, 4.00)( 8.50, 3.00)
  \qbezier(10.00, 0.00)(10.00, 1.00)( 7.50, 2.50)
  \qbezier( 5.00, 5.00)( 5.00, 4.00)( 7.50, 2.50)
  \qbezier( 0.00, 5.00)( 0.00, 6.00)( 2.50, 7.50)
  \qbezier( 5.00,10.00)( 5.00, 9.00)( 2.50, 7.50)
  \qbezier( 5.00, 5.00)( 5.00, 6.00)( 3.50, 7.00)
  \qbezier( 0.00,10.00)( 0.00, 9.00)( 1.50, 8.00)
  \qbezier( 5.00,10.00)( 5.00,11.00)( 6.50,12.00)
  \qbezier(10.00,15.00)(10.00,14.00)( 8.50,13.00)
  \qbezier(10.00,10.00)(10.00,11.00)( 7.50,12.50)
  \qbezier( 5.00,15.00)( 5.00,14.00)( 7.50,12.50)
  \put(10.00, 5.00){\line( 0, 1){ 5.00}}
  \put( 0.00, 0.00){\line( 0, 1){ 5.00}}
  \put(10.00,15.00){\vector( 0, 1){ 2.00}}
  \put( 5.00,15.00){\vector( 0, 1){ 2.00}}
  \put( 0.00,10.00){\vector( 0, 1){ 7.00}}
  \put( 0.00,-0.50){\makebox( 0.00, 0.00)[ct]{#1}}
  \put( 5.00,-0.50){\makebox( 0.00, 0.00)[ct]{#2}}
  \put(10.00,-0.50){\makebox( 0.00, 0.00)[ct]{#3}}
  \end{picture}}}

\newcommand{\ThreeBraidNTwoNOnePTwo}[3]
  {\raisebox{-6.75mm}
  {\begin{picture}(10.00,15.00)
  \qbezier( 5.00, 0.00)( 5.00, 1.00)( 6.50, 2.00)
  \qbezier(10.00, 5.00)(10.00, 4.00)( 8.50, 3.00)
  \qbezier(10.00, 0.00)(10.00, 1.00)( 7.50, 2.50)
  \qbezier( 5.00, 5.00)( 5.00, 4.00)( 7.50, 2.50)
  \qbezier( 0.00, 5.00)( 0.00, 6.00)( 1.50, 7.00)
  \qbezier( 5.00,10.00)( 5.00, 9.00)( 3.50, 8.00)
  \qbezier( 5.00, 5.00)( 5.00, 6.00)( 2.50, 7.50)
  \qbezier( 0.00,10.00)( 0.00, 9.00)( 2.50, 7.50)
  \qbezier( 5.00,10.00)( 5.00,11.00)( 7.50,12.50)
  \qbezier(10.00,15.00)(10.00,14.00)( 7.50,12.50)
  \qbezier(10.00,10.00)(10.00,11.00)( 8.50,12.00)
  \qbezier( 5.00,15.00)( 5.00,14.00)( 6.50,13.00)
  \put(10.00, 5.00){\line( 0, 1){ 5.00}}
  \put( 0.00, 0.00){\line( 0, 1){ 5.00}}
  \put(10.00,15.00){\vector( 0, 1){ 2.00}}
  \put( 5.00,15.00){\vector( 0, 1){ 2.00}}
  \put( 0.00,10.00){\vector( 0, 1){ 7.00}}
  \put( 0.00,-0.50){\makebox( 0.00, 0.00)[ct]{#1}}
  \put( 5.00,-0.50){\makebox( 0.00, 0.00)[ct]{#2}}
  \put(10.00,-0.50){\makebox( 0.00, 0.00)[ct]{#3}}
  \end{picture}}}

\newcommand{\ThreeBraidNTwoNOneNTwo}[3]
  {\raisebox{-6.75mm}
  {\begin{picture}(10.00,15.00)
  \qbezier( 5.00, 0.00)( 5.00, 1.00)( 6.50, 2.00)
  \qbezier(10.00, 5.00)(10.00, 4.00)( 8.50, 3.00)
  \qbezier(10.00, 0.00)(10.00, 1.00)( 7.50, 2.50)
  \qbezier( 5.00, 5.00)( 5.00, 4.00)( 7.50, 2.50)
  \qbezier( 0.00, 5.00)( 0.00, 6.00)( 1.50, 7.00)
  \qbezier( 5.00,10.00)( 5.00, 9.00)( 3.50, 8.00)
  \qbezier( 5.00, 5.00)( 5.00, 6.00)( 2.50, 7.50)
  \qbezier( 0.00,10.00)( 0.00, 9.00)( 2.50, 7.50)
  \qbezier( 5.00,10.00)( 5.00,11.00)( 6.50,12.00)
  \qbezier(10.00,15.00)(10.00,14.00)( 8.50,13.00)
  \qbezier(10.00,10.00)(10.00,11.00)( 7.50,12.50)
  \qbezier( 5.00,15.00)( 5.00,14.00)( 7.50,12.50)
  \put(10.00, 5.00){\line( 0, 1){ 5.00}}
  \put( 0.00, 0.00){\line( 0, 1){ 5.00}}
  \put(10.00,15.00){\vector( 0, 1){ 2.00}}
  \put( 5.00,15.00){\vector( 0, 1){ 2.00}}
  \put( 0.00,10.00){\vector( 0, 1){ 7.00}}
  \put( 0.00,-0.50){\makebox( 0.00, 0.00)[ct]{#1}}
  \put( 5.00,-0.50){\makebox( 0.00, 0.00)[ct]{#2}}
  \put(10.00,-0.50){\makebox( 0.00, 0.00)[ct]{#3}}
  \end{picture}}}

\newcommand{\ThreeBraidCOneCTwoCOne}[3]
  {\raisebox{-6.75mm}
  {\begin{picture}(10.00,15.00)
  \qbezier( 0.00, 0.00)( 0.00, 1.00)( 2.50, 2.50)
  \qbezier( 5.00, 5.00)( 5.00, 4.00)( 2.50, 2.50)
  \qbezier( 5.00, 0.00)( 5.00, 1.00)( 2.50, 2.50)
  \qbezier( 0.00, 5.00)( 0.00, 4.00)( 2.50, 2.50)
  \qbezier( 5.00, 5.00)( 5.00, 6.00)( 7.50, 7.50)
  \qbezier(10.00,10.00)(10.00, 9.00)( 7.50, 7.50)
  \qbezier(10.00, 5.00)(10.00, 6.00)( 7.50, 7.50)
  \qbezier( 5.00,10.00)( 5.00, 9.00)( 7.50, 7.50)
  \qbezier( 0.00,10.00)( 0.00,11.00)( 2.50,12.50)
  \qbezier( 5.00,15.00)( 5.00,14.00)( 2.50,12.50)
  \qbezier( 5.00,10.00)( 5.00,11.00)( 2.50,12.50)
  \qbezier( 0.00,15.00)( 0.00,14.00)( 2.50,12.50)
  \put( 0.00, 5.00){\line( 0, 1){ 5.00}}
  \put(10.00, 0.00){\line( 0, 1){ 5.00}}
  \put( 0.00,15.00){\vector( 0, 1){ 2.00}}
  \put( 5.00,15.00){\vector( 0, 1){ 2.00}}
  \put(10.00,10.00){\vector( 0, 1){ 7.00}}
  \put( 2.50, 2.50){\circle*{ 2.00}}
  \put( 7.50, 7.50){\circle*{ 2.00}}
  \put( 2.50,12.50){\circle*{ 2.00}}
  \put( 0.00,-0.50){\makebox( 0.00, 0.00)[ct]{#1}}
  \put( 5.00,-0.50){\makebox( 0.00, 0.00)[ct]{#2}}
  \put(10.00,-0.50){\makebox( 0.00, 0.00)[ct]{#3}}
  \end{picture}}}

\newcommand{\ThreeBraidCTwoCOneCTwo}[3]
  {\raisebox{-6.75mm}
  {\begin{picture}(10.00,15.00)
  \qbezier( 5.00, 0.00)( 5.00, 1.00)( 7.50, 2.50)
  \qbezier(10.00, 5.00)(10.00, 4.00)( 7.50, 2.50)
  \qbezier(10.00, 0.00)(10.00, 1.00)( 7.50, 2.50)
  \qbezier( 5.00, 5.00)( 5.00, 4.00)( 7.50, 2.50)
  \qbezier( 0.00, 5.00)( 0.00, 6.00)( 2.50, 7.50)
  \qbezier( 5.00,10.00)( 5.00, 9.00)( 2.50, 7.50)
  \qbezier( 5.00, 5.00)( 5.00, 6.00)( 2.50, 7.50)
  \qbezier( 0.00,10.00)( 0.00, 9.00)( 2.50, 7.50)
  \qbezier( 5.00,10.00)( 5.00,11.00)( 7.50,12.50)
  \qbezier(10.00,15.00)(10.00,14.00)( 7.50,12.50)
  \qbezier(10.00,10.00)(10.00,11.00)( 7.50,12.50)
  \qbezier( 5.00,15.00)( 5.00,14.00)( 7.50,12.50)
  \put( 7.50, 2.50){\circle*{ 2.00}}
  \put( 2.50, 7.50){\circle*{ 2.00}}
  \put( 7.50,12.50){\circle*{ 2.00}}
  \put(10.00, 5.00){\line( 0, 1){ 5.00}}
  \put( 0.00, 0.00){\line( 0, 1){ 5.00}}
  \put(10.00,15.00){\vector( 0, 1){ 2.00}}
  \put( 5.00,15.00){\vector( 0, 1){ 2.00}}
  \put( 0.00,10.00){\vector( 0, 1){ 7.00}}
  \put( 0.00,-0.50){\makebox( 0.00, 0.00)[ct]{#1}}
  \put( 5.00,-0.50){\makebox( 0.00, 0.00)[ct]{#2}}
  \put(10.00,-0.50){\makebox( 0.00, 0.00)[ct]{#3}}
  \end{picture}}}

\newcommand{\ThreeBraidPTwoPOnePOnePTwo}[3]
  {\raisebox{-9.00mm}
  {\begin{picture}(10.00,20.00)
  \qbezier( 5.00, 0.00)( 5.00, 1.00)( 7.50, 2.50)
  \qbezier(10.00, 5.00)(10.00, 4.00)( 7.50, 2.50)
  \qbezier(10.00, 0.00)(10.00, 1.00)( 8.50, 2.00)
  \qbezier( 5.00, 5.00)( 5.00, 4.00)( 6.50, 3.00)
  \qbezier( 0.00, 5.00)( 0.00, 6.00)( 2.50, 7.50)
  \qbezier( 5.00,10.00)( 5.00, 9.00)( 2.50, 7.50)
  \qbezier( 5.00, 5.00)( 5.00, 6.00)( 3.50, 7.00)
  \qbezier( 0.00,10.00)( 0.00, 9.00)( 1.50, 8.00)
  \qbezier( 0.00,10.00)( 0.00,11.00)( 2.50,12.50)
  \qbezier( 5.00,15.00)( 5.00,14.00)( 2.50,12.50)
  \qbezier( 5.00,10.00)( 5.00,11.00)( 3.50,12.00)
  \qbezier( 0.00,15.00)( 0.00,14.00)( 1.50,13.00)
  \qbezier( 5.00,15.00)( 5.00,16.00)( 7.50,17.50)
  \qbezier(10.00,20.00)(10.00,19.00)( 7.50,17.50)
  \qbezier(10.00,15.00)(10.00,16.00)( 8.50,17.00)
  \qbezier( 5.00,20.00)( 5.00,19.00)( 6.50,18.00)
  \put(10.00, 5.00){\line( 0, 1){10.00}}
  \put( 0.00, 0.00){\line( 0, 1){ 5.00}}
  \put(10.00,20.00){\vector( 0, 1){ 2.00}}
  \put( 5.00,20.00){\vector( 0, 1){ 2.00}}
  \put( 0.00,15.00){\vector( 0, 1){ 7.00}}
  \put( 0.00,-0.50){\makebox( 0.00, 0.00)[ct]{#1}}
  \put( 5.00,-0.50){\makebox( 0.00, 0.00)[ct]{#2}}
  \put(10.00,-0.50){\makebox( 0.00, 0.00)[ct]{#3}}
  \end{picture}}}

\newcommand{\ThreeBraidPOnePTwoPTwoPOne}[3]
  {\raisebox{-9.00mm}
  {\begin{picture}(10.00,20.00)
  \qbezier( 0.00, 0.00)( 0.00, 1.00)( 2.50, 2.50)
  \qbezier( 5.00, 5.00)( 5.00, 4.00)( 2.50, 2.50)
  \qbezier( 5.00, 0.00)( 5.00, 1.00)( 3.50, 2.00)
  \qbezier( 0.00, 5.00)( 0.00, 4.00)( 1.50, 3.00)
  \qbezier( 5.00, 5.00)( 5.00, 6.00)( 7.50, 7.50)
  \qbezier(10.00,10.00)(10.00, 9.00)( 7.50, 7.50)
  \qbezier(10.00, 5.00)(10.00, 6.00)( 8.50, 7.00)
  \qbezier( 5.00,10.00)( 5.00, 9.00)( 6.50, 8.00)
  \qbezier( 5.00,10.00)( 5.00,11.00)( 7.50,12.50)
  \qbezier(10.00,15.00)(10.00,14.00)( 7.50,12.50)
  \qbezier(10.00,10.00)(10.00,11.00)( 8.50,12.00)
  \qbezier( 5.00,15.00)( 5.00,14.00)( 6.50,13.00)
  \qbezier( 0.00,15.00)( 0.00,16.00)( 2.50,17.50)
  \qbezier( 5.00,20.00)( 5.00,19.00)( 2.50,17.50)
  \qbezier( 5.00,15.00)( 5.00,16.00)( 3.50,17.00)
  \qbezier( 0.00,20.00)( 0.00,19.00)( 1.50,18.00)
  \put( 0.00, 5.00){\line( 0, 1){10.00}}
  \put(10.00, 0.00){\line( 0, 1){ 5.00}}
  \put( 0.00,20.00){\vector( 0, 1){ 2.00}}
  \put( 5.00,20.00){\vector( 0, 1){ 2.00}}
  \put(10.00,15.00){\vector( 0, 1){ 7.00}}
  \put( 0.00,-0.50){\makebox( 0.00, 0.00)[ct]{#1}}
  \put( 5.00,-0.50){\makebox( 0.00, 0.00)[ct]{#2}}
  \put(10.00,-0.50){\makebox( 0.00, 0.00)[ct]{#3}}
  \end{picture}}}

\newcommand{\ThreeBraidPOnePOnePTwoPTwo}[3]
  {\raisebox{-9.00mm}
  {\begin{picture}(10.00,20.00)
  \qbezier( 0.00, 0.00)( 0.00, 1.00)( 2.50, 2.50)
  \qbezier( 5.00, 5.00)( 5.00, 4.00)( 2.50, 2.50)
  \qbezier( 5.00, 0.00)( 5.00, 1.00)( 3.50, 2.00)
  \qbezier( 0.00, 5.00)( 0.00, 4.00)( 1.50, 3.00)
  \qbezier( 0.00, 5.00)( 0.00, 6.00)( 2.50, 7.50)
  \qbezier( 5.00,10.00)( 5.00, 9.00)( 2.50, 7.50)
  \qbezier( 5.00, 5.00)( 5.00, 6.00)( 3.50, 7.00)
  \qbezier( 0.00,10.00)( 0.00, 9.00)( 1.50, 8.00)
  \qbezier( 5.00,10.00)( 5.00,11.00)( 7.50,12.50)
  \qbezier(10.00,15.00)(10.00,14.00)( 7.50,12.50)
  \qbezier(10.00,10.00)(10.00,11.00)( 8.50,12.00)
  \qbezier( 5.00,15.00)( 5.00,14.00)( 6.50,13.00)
  \qbezier( 5.00,15.00)( 5.00,16.00)( 7.50,17.50)
  \qbezier(10.00,20.00)(10.00,19.00)( 7.50,17.50)
  \qbezier(10.00,15.00)(10.00,16.00)( 8.50,17.00)
  \qbezier( 5.00,20.00)( 5.00,19.00)( 6.50,18.00)
  \put(10.00, 0.00){\line( 0, 1){10.00}}
  \put( 0.00,10.00){\vector( 0, 1){12.00}}
  \put( 5.00,20.00){\vector( 0, 1){ 2.00}}
  \put(10.00,20.00){\vector( 0, 1){ 2.00}}
  \put( 0.00,-0.50){\makebox( 0.00, 0.00)[ct]{#1}}
  \put( 5.00,-0.50){\makebox( 0.00, 0.00)[ct]{#2}}
  \put(10.00,-0.50){\makebox( 0.00, 0.00)[ct]{#3}}
  \end{picture}}}

\newcommand{\ThreeBraidPTwoPTwoPOnePOne}[3]
  {\raisebox{-9.00mm}
  {\begin{picture}(10.00,20.00)
  \qbezier( 5.00, 0.00)( 5.00, 1.00)( 7.50, 2.50)
  \qbezier(10.00, 5.00)(10.00, 4.00)( 7.50, 2.50)
  \qbezier(10.00, 0.00)(10.00, 1.00)( 8.50, 2.00)
  \qbezier( 5.00, 5.00)( 5.00, 4.00)( 6.50, 3.00)
  \qbezier( 5.00, 5.00)( 5.00, 6.00)( 7.50, 7.50)
  \qbezier(10.00,10.00)(10.00, 9.00)( 7.50, 7.50)
  \qbezier(10.00, 5.00)(10.00, 6.00)( 8.50, 7.00)
  \qbezier( 5.00,10.00)( 5.00, 9.00)( 6.50, 8.00)
  \qbezier( 0.00,10.00)( 0.00,11.00)( 2.50,12.50)
  \qbezier( 5.00,15.00)( 5.00,14.00)( 2.50,12.50)
  \qbezier( 5.00,10.00)( 5.00,11.00)( 3.50,12.00)
  \qbezier( 0.00,15.00)( 0.00,14.00)( 1.50,13.00)
  \qbezier( 0.00,15.00)( 0.00,16.00)( 2.50,17.50)
  \qbezier( 5.00,20.00)( 5.00,19.00)( 2.50,17.50)
  \qbezier( 5.00,15.00)( 5.00,16.00)( 3.50,17.00)
  \qbezier( 0.00,20.00)( 0.00,19.00)( 1.50,18.00)
  \put( 0.00, 0.00){\line( 0, 1){10.00}}
  \put(10.00,10.00){\vector( 0, 1){12.00}}
  \put( 5.00,20.00){\vector( 0, 1){ 2.00}}
  \put( 0.00,20.00){\vector( 0, 1){ 2.00}}
  \put( 0.00,-0.50){\makebox( 0.00, 0.00)[ct]{#1}}
  \put( 5.00,-0.50){\makebox( 0.00, 0.00)[ct]{#2}}
  \put(10.00,-0.50){\makebox( 0.00, 0.00)[ct]{#3}}
  \end{picture}}}

\newcommand{\ThreeBraidPOnePOne}[3]
  {\raisebox{-9.00mm}
  {\begin{picture}(10.00,20.00)
  \qbezier( 0.00, 5.00)( 0.00, 6.00)( 2.50, 7.50)
  \qbezier( 5.00,10.00)( 5.00, 9.00)( 2.50, 7.50)
  \qbezier( 5.00, 5.00)( 5.00, 6.00)( 3.50, 7.00)
  \qbezier( 0.00,10.00)( 0.00, 9.00)( 1.50, 8.00)
  \qbezier( 0.00,10.00)( 0.00,11.00)( 2.50,12.50)
  \qbezier( 5.00,15.00)( 5.00,14.00)( 2.50,12.50)
  \qbezier( 5.00,10.00)( 5.00,11.00)( 3.50,12.00)
  \qbezier( 0.00,15.00)( 0.00,14.00)( 1.50,13.00)
  \put( 0.00, 0.00){\line( 0, 1){ 5.00}}
  \put( 5.00, 0.00){\line( 0, 1){ 5.00}}
  \put( 0.00,15.00){\vector( 0, 1){ 7.00}}
  \put( 5.00,15.00){\vector( 0, 1){ 7.00}}
  \put(10.00, 0.00){\vector( 0, 1){22.00}}
  \put( 0.00,-0.50){\makebox( 0.00, 0.00)[ct]{#1}}
  \put( 5.00,-0.50){\makebox( 0.00, 0.00)[ct]{#2}}
  \put(10.00,-0.50){\makebox( 0.00, 0.00)[ct]{#3}}
  \end{picture}}}

\newcommand{\ThreeBraidPTwoPTwo}[3]
  {\raisebox{-9.00mm}
  {\begin{picture}(10.00,20.00)
  \qbezier( 5.00, 5.00)( 5.00, 6.00)( 7.50, 7.50)
  \qbezier(10.00,10.00)(10.00, 9.00)( 7.50, 7.50)
  \qbezier(10.00, 5.00)(10.00, 6.00)( 8.50, 7.00)
  \qbezier( 5.00,10.00)( 5.00, 9.00)( 6.50, 8.00)
  \qbezier( 5.00,10.00)( 5.00,11.00)( 7.50,12.50)
  \qbezier(10.00,15.00)(10.00,14.00)( 7.50,12.50)
  \qbezier(10.00,10.00)(10.00,11.00)( 8.50,12.00)
  \qbezier( 5.00,15.00)( 5.00,14.00)( 6.50,13.00)
  \put(10.00, 0.00){\line( 0, 1){ 5.00}}
  \put( 5.00, 0.00){\line( 0, 1){ 5.00}}
  \put(10.00,15.00){\vector( 0, 1){ 7.00}}
  \put( 5.00,15.00){\vector( 0, 1){ 7.00}}
  \put( 0.00, 0.00){\vector( 0, 1){22.00}}
  \put( 0.00,-0.50){\makebox( 0.00, 0.00)[ct]{#1}}
  \put( 5.00,-0.50){\makebox( 0.00, 0.00)[ct]{#2}}
  \put(10.00,-0.50){\makebox( 0.00, 0.00)[ct]{#3}}
  \end{picture}}}

\newcommand{\ThreeBraidPTwoPOneCOneCTwo}[3]
  {\raisebox{-9.00mm}
  {\begin{picture}(10.00,20.00)
  \qbezier( 5.00, 0.00)( 5.00, 1.00)( 7.50, 2.50)
  \qbezier(10.00, 5.00)(10.00, 4.00)( 7.50, 2.50)
  \qbezier(10.00, 0.00)(10.00, 1.00)( 8.50, 2.00)
  \qbezier( 5.00, 5.00)( 5.00, 4.00)( 6.50, 3.00)
  \qbezier( 0.00, 5.00)( 0.00, 6.00)( 2.50, 7.50)
  \qbezier( 5.00,10.00)( 5.00, 9.00)( 2.50, 7.50)
  \qbezier( 5.00, 5.00)( 5.00, 6.00)( 3.50, 7.00)
  \qbezier( 0.00,10.00)( 0.00, 9.00)( 1.50, 8.00)
  \qbezier( 0.00,10.00)( 0.00,11.00)( 2.50,12.50)
  \qbezier( 5.00,15.00)( 5.00,14.00)( 2.50,12.50)
  \qbezier( 5.00,10.00)( 5.00,11.00)( 2.50,12.50)
  \qbezier( 0.00,15.00)( 0.00,14.00)( 2.50,12.50)
  \qbezier( 5.00,15.00)( 5.00,16.00)( 7.50,17.50)
  \qbezier(10.00,20.00)(10.00,19.00)( 7.50,17.50)
  \qbezier(10.00,15.00)(10.00,16.00)( 7.50,17.50)
  \qbezier( 5.00,20.00)( 5.00,19.00)( 7.50,17.50)
  \put( 2.50,12.50){\circle*{ 2.00}}
  \put( 7.50,17.50){\circle*{ 2.00}}
  \put(10.00, 5.00){\line( 0, 1){10.00}}
  \put( 0.00, 0.00){\line( 0, 1){ 5.00}}
  \put(10.00,20.00){\vector( 0, 1){ 2.00}}
  \put( 5.00,20.00){\vector( 0, 1){ 2.00}}
  \put( 0.00,15.00){\vector( 0, 1){ 7.00}}
  \put( 0.00,-0.50){\makebox( 0.00, 0.00)[ct]{#1}}
  \put( 5.00,-0.50){\makebox( 0.00, 0.00)[ct]{#2}}
  \put(10.00,-0.50){\makebox( 0.00, 0.00)[ct]{#3}}
  \end{picture}}}

\newcommand{\ThreeBraidCOneCTwoPTwoPOne}[3]
  {\raisebox{-9.00mm}
  {\begin{picture}(10.00,20.00)
  \qbezier( 0.00, 0.00)( 0.00, 1.00)( 2.50, 2.50)
  \qbezier( 5.00, 5.00)( 5.00, 4.00)( 2.50, 2.50)
  \qbezier( 5.00, 0.00)( 5.00, 1.00)( 2.50, 2.50)
  \qbezier( 0.00, 5.00)( 0.00, 4.00)( 2.50, 2.50)
  \qbezier( 5.00, 5.00)( 5.00, 6.00)( 7.50, 7.50)
  \qbezier(10.00,10.00)(10.00, 9.00)( 7.50, 7.50)
  \qbezier(10.00, 5.00)(10.00, 6.00)( 7.50, 7.50)
  \qbezier( 5.00,10.00)( 5.00, 9.00)( 7.50, 7.50)
  \qbezier( 5.00,10.00)( 5.00,11.00)( 7.50,12.50)
  \qbezier(10.00,15.00)(10.00,14.00)( 7.50,12.50)
  \qbezier(10.00,10.00)(10.00,11.00)( 8.50,12.00)
  \qbezier( 5.00,15.00)( 5.00,14.00)( 6.50,13.00)
  \qbezier( 0.00,15.00)( 0.00,16.00)( 2.50,17.50)
  \qbezier( 5.00,20.00)( 5.00,19.00)( 2.50,17.50)
  \qbezier( 5.00,15.00)( 5.00,16.00)( 3.50,17.00)
  \qbezier( 0.00,20.00)( 0.00,19.00)( 1.50,18.00)
  \put( 2.50, 2.50){\circle*{ 2.00}}
  \put( 7.50, 7.50){\circle*{ 2.00}}
  \put( 0.00, 5.00){\line( 0, 1){10.00}}
  \put(10.00, 0.00){\line( 0, 1){ 5.00}}
  \put( 0.00,20.00){\vector( 0, 1){ 2.00}}
  \put( 5.00,20.00){\vector( 0, 1){ 2.00}}
  \put(10.00,15.00){\vector( 0, 1){ 7.00}}
  \put( 0.00,-0.50){\makebox( 0.00, 0.00)[ct]{#1}}
  \put( 5.00,-0.50){\makebox( 0.00, 0.00)[ct]{#2}}
  \put(10.00,-0.50){\makebox( 0.00, 0.00)[ct]{#3}}
  \end{picture}}}

\newcommand{\ThreeBraidPOneCOnePTwoCTwo}[3]
  {\raisebox{-9.00mm}
  {\begin{picture}(10.00,20.00)
  \qbezier( 0.00, 0.00)( 0.00, 1.00)( 2.50, 2.50)
  \qbezier( 5.00, 5.00)( 5.00, 4.00)( 2.50, 2.50)
  \qbezier( 5.00, 0.00)( 5.00, 1.00)( 3.50, 2.00)
  \qbezier( 0.00, 5.00)( 0.00, 4.00)( 1.50, 3.00)
  \qbezier( 0.00, 5.00)( 0.00, 6.00)( 2.50, 7.50)
  \qbezier( 5.00,10.00)( 5.00, 9.00)( 2.50, 7.50)
  \qbezier( 5.00, 5.00)( 5.00, 6.00)( 2.50, 7.50)
  \qbezier( 0.00,10.00)( 0.00, 9.00)( 2.50, 7.50)
  \qbezier( 5.00,10.00)( 5.00,11.00)( 7.50,12.50)
  \qbezier(10.00,15.00)(10.00,14.00)( 7.50,12.50)
  \qbezier(10.00,10.00)(10.00,11.00)( 8.50,12.00)
  \qbezier( 5.00,15.00)( 5.00,14.00)( 6.50,13.00)
  \qbezier( 5.00,15.00)( 5.00,16.00)( 7.50,17.50)
  \qbezier(10.00,20.00)(10.00,19.00)( 7.50,17.50)
  \qbezier(10.00,15.00)(10.00,16.00)( 7.50,17.50)
  \qbezier( 5.00,20.00)( 5.00,19.00)( 7.50,17.50)
  \put( 2.50, 7.50){\circle*{ 2.00}}
  \put( 7.50,17.50){\circle*{ 2.00}}
  \put(10.00, 0.00){\line( 0, 1){10.00}}
  \put( 0.00,10.00){\vector( 0, 1){12.00}}
  \put( 5.00,20.00){\vector( 0, 1){ 2.00}}
  \put(10.00,20.00){\vector( 0, 1){ 2.00}}
  \put( 0.00,-0.50){\makebox( 0.00, 0.00)[ct]{#1}}
  \put( 5.00,-0.50){\makebox( 0.00, 0.00)[ct]{#2}}
  \put(10.00,-0.50){\makebox( 0.00, 0.00)[ct]{#3}}
  \end{picture}}}

\newcommand{\ThreeBraidPTwoCTwoPOneCOne}[3]
  {\raisebox{-9.00mm}
  {\begin{picture}(10.00,20.00)
  \qbezier( 5.00, 0.00)( 5.00, 1.00)( 7.50, 2.50)
  \qbezier(10.00, 5.00)(10.00, 4.00)( 7.50, 2.50)
  \qbezier(10.00, 0.00)(10.00, 1.00)( 8.50, 2.00)
  \qbezier( 5.00, 5.00)( 5.00, 4.00)( 6.50, 3.00)
  \qbezier( 5.00, 5.00)( 5.00, 6.00)( 7.50, 7.50)
  \qbezier(10.00,10.00)(10.00, 9.00)( 7.50, 7.50)
  \qbezier(10.00, 5.00)(10.00, 6.00)( 7.50, 7.50)
  \qbezier( 5.00,10.00)( 5.00, 9.00)( 7.50, 7.50)
  \qbezier( 0.00,10.00)( 0.00,11.00)( 2.50,12.50)
  \qbezier( 5.00,15.00)( 5.00,14.00)( 2.50,12.50)
  \qbezier( 5.00,10.00)( 5.00,11.00)( 3.50,12.00)
  \qbezier( 0.00,15.00)( 0.00,14.00)( 1.50,13.00)
  \qbezier( 0.00,15.00)( 0.00,16.00)( 2.50,17.50)
  \qbezier( 5.00,20.00)( 5.00,19.00)( 2.50,17.50)
  \qbezier( 5.00,15.00)( 5.00,16.00)( 2.50,17.50)
  \qbezier( 0.00,20.00)( 0.00,19.00)( 2.50,17.50)
  \put( 7.50, 7.50){\circle*{ 2.00}}
  \put( 2.50,17.50){\circle*{ 2.00}}
  \put( 0.00, 0.00){\line( 0, 1){10.00}}
  \put(10.00,10.00){\vector( 0, 1){12.00}}
  \put( 5.00,20.00){\vector( 0, 1){ 2.00}}
  \put( 0.00,20.00){\vector( 0, 1){ 2.00}}
  \put( 0.00,-0.50){\makebox( 0.00, 0.00)[ct]{#1}}
  \put( 5.00,-0.50){\makebox( 0.00, 0.00)[ct]{#2}}
  \put(10.00,-0.50){\makebox( 0.00, 0.00)[ct]{#3}}
  \end{picture}}}

\newcommand{\PositiveClasp}[2]
  {\raisebox{-6.75mm}
  {\begin{picture}(10.00,15.00)
  \put( 4.00, 7.50){\oval( 8.00, 7.00)[l]}
  \put( 6.00, 7.50){\oval( 8.00, 7.00)[r]}
  \put( 4.00, 4.00){\line( 1, 0){ 2.00}} 
  \put( 5.00, 5.00){\vector( 0, 1){11.00}}
  \put(10.00, 9.00){\vector( 0, 1){ 0.00}}
  \put( 5.00, 0.00){\line( 0, 1){ 3.00}}
  \put( 5.00,-0.50){\makebox( 0.00, 0.00)[ct]{#1}}
  \put(11.00, 7.50){\makebox( 0.00, 0.00)[lc]{#2}}
  \end{picture}}}

\newcommand{\PositiveClaspWithDoublePoint}[2]
  {\raisebox{-6.75mm}
  {\begin{picture}(10.00,15.00)
  \put( 4.00, 7.50){\oval( 8.00, 7.00)[l]}
  \put( 6.00, 7.50){\oval( 8.00, 7.00)[r]}
  \put( 4.00, 4.00){\line( 1, 0){ 2.00}} 
  \put( 5.00, 0.00){\vector( 0, 1){16.00}}
  \put(10.00, 9.00){\vector( 0, 1){ 0.00}}
  \put( 5.00, 4.00){\circle*{ 2.00}}
  \put( 5.00,-0.50){\makebox( 0.00, 0.00)[ct]{#1}}
  \put(11.00, 7.50){\makebox( 0.00, 0.00)[lc]{#2}}
  \end{picture}}}

\newcommand{\WirtingerPositive}
  {\raisebox{-6.75mm}
  {\begin{picture}(15.00,15.00)
  \put( 0.00, 0.00){\vector( 1, 1){15.00}}
  \put(15.00, 0.00){\line(-1, 1){ 6.50}}
  \put( 6.50, 8.50){\vector(-1, 1){ 6.50}}
  \thinlines
  \put( 7.50, 7.50){\circle{10.00}}
  \thicklines
  \put(12.50, 7.50){\circle*{ 1.00}}
  \put(11.60,10.50){\vector(-1, 3){ 0.00}}
  \put(13.50, 7.50){\makebox( 0.00, 0.00)[lc]{$p_i$}}
  \put( 7.50, 6.00){\makebox( 0.00, 0.00)[ct]{$c_i$}}
  \put( 0.00,15.00){\makebox( 0.00, 0.00)[rb]{$x_i$}}
  \put(15.00,15.00){\makebox( 0.00, 0.00)[lb]{$x_j$}}
  \put(15.00, 0.00){\makebox( 0.00, 0.00)[lt]{$x_k$}}
  \end{picture}}}
  
\newcommand{\WirtingerNegative}
  {\raisebox{-6.75mm}
  {\begin{picture}(15.00,15.00)
  \put(15.00, 0.00){\vector(-1, 1){15.00}}
  \put( 0.00, 0.00){\line( 1, 1){ 6.50}}
  \put( 8.50, 8.50){\vector( 1, 1){ 6.50}}
  \thinlines
  \put( 7.50, 7.50){\circle{10.00}}
  \thicklines
  \put(12.50, 7.50){\circle*{ 1.00}}
  \put(11.60,10.50){\vector(-1, 3){ 0.00}}
  \put(13.50, 7.50){\makebox( 0.00, 0.00)[lc]{$p_i$}}
  \put( 7.50, 6.00){\makebox( 0.00, 0.00)[ct]{$c_i$}}
  \put( 0.00,15.00){\makebox( 0.00, 0.00)[rb]{$x_j$}}
  \put(15.00,15.00){\makebox( 0.00, 0.00)[lb]{$x_i$}}
  \put( 0.00, 0.00){\makebox( 0.00, 0.00)[rt]{$x_k$}}
  \end{picture}}}
  
\newcommand{\EllPlus}
  {\raisebox{-6.75mm}
  {\begin{picture}(15.00,15.00)
  \put( 0.00,15.00){\line( 1,-1){ 4.00}}
  \put( 5.20, 9.80){\vector( 1,-1){ 9.80}}
  \put( 0.00, 0.00){\line( 1, 1){ 6.50}}
  \put( 8.50, 8.50){\vector( 1, 1){ 6.50}}
  \qbezier( 4.00,11.00)( 5.00,10.00)( 6.00,11.00)
  \qbezier( 6.00,11.00)( 7.00,12.00)( 6.00,13.00)
  \qbezier( 6.00,13.00)( 5.00,14.00)( 4.00,13.00)
  \qbezier( 4.00,13.00)( 3.70,12.70)( 4.00,12.00)
  \put( 7.50, 6.00){\makebox( 0.00, 0.00)[ct]{$c_2$}}
  \put( 3.80,10.80){\makebox( 0.00, 0.00)[rt]{$c_1$}}
  \put( 0.00, 0.00){\makebox( 0.00, 0.00)[rt]{$x_4$}}
  \put( 0.00,15.00){\makebox( 0.00, 0.00)[rb]{$x_3$}}
  \put(15.00,15.00){\makebox( 0.00, 0.00)[lb]{$x_2$}}
  \put(15.00, 0.00){\makebox( 0.00, 0.00)[lt]{$x_1$}}
  \end{picture}}}

\newcommand{\EllMinus}
  {\raisebox{-6.75mm}
  {\begin{picture}(15.00,15.00)
  \put( 0.00, 0.00){\line( 1, 1){ 4.00}}
  \put( 5.20, 5.20){\vector( 1, 1){ 9.80}}
  \put( 0.00,15.00){\line( 1,-1){ 6.50}}
  \put( 8.50, 6.50){\vector( 1,-1){ 6.50}}
  \qbezier( 4.00, 4.00)( 5.00, 5.00)( 4.00, 6.00)
  \qbezier( 4.00, 6.00)( 3.00, 7.00)( 2.00, 6.00)
  \qbezier( 2.00, 6.00)( 1.00, 5.00)( 2.00, 4.00)
  \qbezier( 2.00, 4.00)( 2.30, 3.70)( 3.00, 4.00)
  \put( 9.00, 7.50){\makebox( 0.00, 0.00)[lc]{$c_1$}}
  \put( 4.20, 3.80){\makebox( 0.00, 0.00)[lt]{$c_2$}}
  \put( 0.00, 0.00){\makebox( 0.00, 0.00)[rt]{$x_4$}}
  \put( 0.00,15.00){\makebox( 0.00, 0.00)[rb]{$x_3$}}
  \put(15.00,15.00){\makebox( 0.00, 0.00)[lb]{$x_2$}}
  \put(15.00, 0.00){\makebox( 0.00, 0.00)[lt]{$x_1$}}
  \end{picture}}}

\newcommand{\EllDot}
  {\raisebox{-6.75mm}
  {\begin{picture}(15.00,15.00)
  \put( 0.00, 0.00){\vector( 1, 1){15.00}}
  \put( 0.00,15.00){\vector( 1,-1){15.00}}
  \put( 7.50, 7.50){\circle*{ 2.00}}
  \put( 0.00, 0.00){\makebox( 0.00, 0.00)[rt]{$x_4$}}
  \put( 0.00,15.00){\makebox( 0.00, 0.00)[rb]{$x_3$}}
  \put(15.00,15.00){\makebox( 0.00, 0.00)[lb]{$x_2$}}
  \put(15.00, 0.00){\makebox( 0.00, 0.00)[lt]{$x_1$}}
  \end{picture}}}

\newcommand{\OneTerm}
  {\raisebox{-6.75mm}
  {\begin{picture}(10.00,10.00)
  \put( 5.00, 0.00){\vector( 0, 1){16.00}}
  \qbezier[10]( 5.00, 5.00)(10.00, 7.50)( 5.00,10.00)
  \end{picture}}}

\section{Introduction}
In \cite{BarNatan/Garoufalidis:MMR_conjecture} D.~Bar-Natan and
S.~Garoufalidis used a weight system for the Alexander
polynomial of knots to prove the so-called Melvin--Morton--Rozansky conjecture
\cite{Melvin/Morton:coloured_Jones,Rozansky:trivial_connection_I}
which relates the Alexander polynomial and some coefficients of coloured Jones
polynomial.
Their weight system can be easily extended to the case of links.
It is a natural problem to give a weight system for the multivariable
Alexander polynomial or its normalised version, the Conway potential function
for links.
\par
The Conway potential function was first introduced by J.H.~Conway
\cite{Conway:enumeration} by giving some `axioms'.
Unfortunately, his `axioms' are not sufficient to define his potential
function.
R.~Hartley \cite{Hartley:Conway_polynomial} gave its precise definition 
by using R.H.~Fox's free differential calculus
\cite{Fox:free_differential_calculus_I,Fox:free_differential_calculus_II}.
He also showed that for two-bridge links Conway's first two identities and
initial data for the trivial knot, for split links, and for the positive Hopf
link are sufficient to calculate the potential function.
After that M.E.~Kidwell \cite{Kidwell:two-variable_Conway} proved they are
also sufficient for calculation of links with two labels $K=T\cup L$ where
$T$ is an unknotted circle and $T$ and $L$ have different labels.
(Note that $L$ may have more than one component.)
Then Y.~Nakanishi \cite{Nakanishi:three-variable_Conway} proved
that we can calculate the potential function if the number of labels
(which equals the number of variables) is two or three.
Besides Hartley's axioms we need Conway's third identity and initial data
for the connected sum of two positive Hopf links and for the three-component
positive Hopf link.
Finally J.~Murakami proved that Conway's first and second identities,
a connect sum formula for the positive Hopf link, initial data for the
trivial knot and for split links, and his new relation involving seven
locally distinct links are sufficient to calculate the Conway potential
function for links with any number of labels.
\par
In this paper we use J.~Murakami's result to define a weight system.
Moreover we will show that our weight system can be calculated recursively
by using five axioms.
Since the proof is fairly easy we expect that there may be similar weight
systems.
If so, by using M.~Kontsevich's integral
\cite{Kontsevich:Vassiliev_invariant} we could then define invariants for
labelled links other than the Conway potential function.
It is also an interesting problem whether our weight system is canonical,
that is, whether we obtain the Conway potential function again after applying
the composition of the Kontsevich integral and our weight system to links.
\begin{acknowledgements}
Most of the results here were prepared for the Knot Seminar at the University
of Liverpool in the autumn term, 1995.
The author thanks H.R.~Morton for inviting him to Liverpool as an EPSRC
fellow, for arranging the Seminar, and for valuable comments.
He also thanks the audience of the Seminar for their patience with his
poor {\em Engrish}.
Thanks are also due to Y.~Nakanishi for helpful conversations.
\end{acknowledgements}
\section{Preliminaries}
In this section we define the Conway potential function and describe the
concept of Vassiliev invariants.
\par
In \cite{Conway:enumeration} J.H.~Conway introduced a notion of
the potential function by an axiomatic way.
R.~Hartley \cite{Hartley:Conway_polynomial} gave its precise definition by
using R.H.~Fox's free differential calculus
\cite{Fox:free_differential_calculus_I,Fox:free_differential_calculus_II} and
proved its existence explicitly.
We follow R.~Hartley to give the definition of the Conway potential function.
\par
Let
$L=K_{1}\cup K_{2}\cup\dots\cup K_{\mu}$ be an oriented $\mu$-component link
with labels $1,2,\dots,n$.
Here $K_{i}$ is labeled with $\ell(i)\in\{1,2,\dots,n\}$ ($i=1,2,\dots,\mu$).  
Let $\Ell$ be a {\em connected} link diagram of $L$.
Let $c_1,c_2,\dots,c_m$ be the crossings and $x_i$ the arc starting at $c_i$
($i=1,2,\dots,m$).
We read the Wirtiger relation $r_i$ at $c_i$ anticlockwise starting at
a point $p_i$ to the right of both arcs.
Then it is of the form
$r_i=x_j x_i x_j^{-1} x_k^{-1}$ or $r_i=x_i x_j x_k^{-1} x_j^{-1}$
according as $c_i$ is a positive crossing or a negative one.
See Figure~1.
\vskip 5mm
\begin{gather*}
  \begin{matrix}
    \WirtingerPositive \\ \\[5mm] \mbox{positive crossing}
  \end{matrix}
  \qquad
  \begin{matrix}
    \WirtingerNegative \\ \\[5mm] \mbox{negative crossing}
  \end{matrix}
  \\[5mm]
  \mbox{Figure 1}
\end{gather*}
\vskip 5mm
\par
Let $\varphi:\ZZ F(x_1,x_2,\dots,x_m)\rightarrow
             \ZZ[{t_1}^{\pm1},{t_2}^{\pm1},\dots,{t_n}^{\pm1}]$ be the
abelianisation homomorphism sending $x_i$ to $t_{\ell(j)}$ if $x_i$
belongs to $K_j$, where $\ZZ F(x_1,x_2,\dots,x_m)$ is the group ring of the
free group generated by $m$ letters $x_1,x_2,\dots,x_m$. 
We consider the $m\times m$ Jacobian matrix
$M(\Ell)=\varphi(\partial r_i/\partial x_j)$,
where $\partial/\partial x_j$ is Fox's free differential calculus
\cite{Fox:free_differential_calculus_I,Fox:free_differential_calculus_II}.
Let $D^{(ij)}(\Ell)$ be the determinant of the matrix obtained from $M(\Ell)$
by deleting the $i$-th row and the $j$-th column.
Then $D^{(ij)}(\Ell)/(\varphi(x_j)-1)$ is the multivariable Alexander
polynomial of the labelled link $L$ (if $\mu>1$; we do not need to
divide by $\varphi(x_j)-1$ if $\mu=1$) and so defines the Conway
potential function up to units of
$\ZZ[{t_1}^{\pm1},{t_2}^{\pm1},\dots,{t_n}^{\pm1}]$.
\par
To define the Conway potential function precisely we need more definitions.
Let $w_i$ be the word read from a path connecting a point in the
unbounded region in $\RR^2$ and $p_i$.
Let $\kappa_g(\Ell)$ be the rotation number (or curvature) of the sublink
consisting of all the components labelled with $g$, which counts
(algebraically) how many times the sublink rotates anticlockwise.
Let $\nu_g(\Ell)$ be the (geometric) number of the crossings where components
labelled with $g$ cross over.
Then according to R.~Hartley we can define the Conway potential function
$\nabla_n(L;t_1,t_2,\dots,t_n)$ as
\begin{equation*}
\nabla_n(L;t_1,t_2,\dots,t_n)
=
\frac{(-1)^{i+j}D^{(ij)}(\Ell)}{\varphi(w_i)(\varphi(x_j)-1)}
\prod_{g=1}^{n}{t_g}^{(\kappa_g(\Ell)-\nu_g(\Ell))/2}.
\end{equation*}
Note that our definition differs from Conway's and Hartley's.
Their potential function is $\nabla_n(L;{t_1}^2,{t_2}^2,\dots,{t_n}^2)$ in
our notation.
\par
It is well known \cite{Fox:free_differential_calculus_II} that
$\nabla_n(L;t_1,t_2,\dots,t_n)\in
\ZZ[{t_1}^{\pm1/2},{t_2}^{\pm1/2},\dots,\linebreak[0]{t_n}^{\pm1/2}]$
if $\mu>1$ and
$({t_1}^{1/2}-{t_1}^{-1/2})\nabla_1(L;t_1)\in\ZZ[{t_1}^{\pm1/2}]$
if $\mu=1$.
In this paper, we study the Laurent expansion of
$\nabla_n(L;\exp(h_1),\exp(h_2),\dots,\linebreak[0]\exp(h_n))$
at 
$(h_1,h_2,\dots,h_n)=(0,0,\dots,0)$ and denote it simply by $\nabla_n(L)$.
So if $\mu>1$, then $\nabla_n(L)$ is a Taylor series and if $\mu=1$, then
it is a Laurent series of the form $\sum_{k=-1}^{\infty}c_k{h_1}^k$.
\par
Next we describe Vassiliev invariants.
Given a numerical link invariant $v$, we can also regard it as an invariant
for singular links, links with double points, as follows.
\begin{equation*}
v\left(\DoublePointUpUp{}{}\right)=
v\left(\PositiveCrossUpUp{}{}\right)-v\left(\NegativeCrossUpUp{}{}\right).
\end{equation*}
Now $v$ is called a Vassiliev invariant of type $d$ 
\cite{Vassiliev:cohomology,Birman/Lin:Vassiliev,Birman:new_point}
if it vanishes for all the singular links with more than $d$ double points.
This is equivalent to saying that for any singular link $L^d$ with $d$ double
points, $v(L^d)$ does not depend on its embedding; it depends only on the
configuration how double points are paired on the circles
\cite{Vassiliev:cohomology,Birman/Lin:Vassiliev,Birman:new_point}.
Such a configuration is described by a chord diagram.
\par
For a compact one-manifold $N$, $N\cup I_1\cup I_2\cup\dots\cup I_m$
is called {\em a chord diagram with support $N$},
where $I_i$ is an interval $[0,1]$ ($1\le i\le m$),
$I_i\cap I_j=\emptyset$ ($i\ne j$),
and $N\cap I_i=N\cap \partial I_i=\partial I_i$ ($1\le i\le m$).
We call a (part of) connected component of $N$ an arc and $I_i$ a chord.
We use solid lines for arcs and dotted lines for chords.
We denote by $\D(N)$ the set of linear combinations of chord diagrams
with support $N$ over $\CC$.
We denote $\D(N)$ modulo the following 4-term relation and the framing
independence relation by $\A(N)$.
\par
\medskip
(4-term relation)
\begin{equation*}
\AmidaAC{}{}{}\,-\,\AmidaCA{}{}{}\,=\,\AmidaCB{}{}{}\,-\,\AmidaBC{}{}{}\,.
\end{equation*}
\par
(framing independence relation)
\begin{equation*}
\OneTerm=0.
\end{equation*}
\par
M.~Kontsevich defined by using the iterated integral a map $Z$ which sends
an embedded circle in $\RR^3$ to an element in $\A(S^1)$
\cite{Kontsevich:Vassiliev_invariant,BarNatan:Vassiliev}.
It is naturally extended to a map sending embedded circles to an element in
$\A(\SplitUnion S^1)$, which is also denoted by $Z$.
His main result is that $Z$ is a link invariant, i.e., it is invariant under
ambient isotopy of $\RR^3$.
So if we have a map $W$ from $\A(\SplitUnion S^1)$ to a ring $R$ then
$W\circ Z$ gives an $R$-valued link invariant.
We call such a $W$ a weight system.
\section{The Conway potential function as Vassiliev invariants}
In this section we will show that every coefficient of $\nabla_n(L)$ is
a Vassiliev invariant.
\begin{definition}
For a Laurent series $f$ with variables $h_1,h_2,\dots,h_n$, we denote by
$c_{p_1,p_2,\dots,p_n}(f)$ the coefficient of 
$\prod_{i=1}^{n}{h_i}^{p_i}$ in $f$.
We also denote by $C_{p}(f)$ the total degree $p$ part of
$f$, which is equal to
$\sum_{\sum p_i=p}\linebreak[0]
c_{p_1,p_2,\dots,p_n}(f)\linebreak[0]
\prod_{i=1}^{n}{h_i}^{p_i}$.
\end{definition}
Then we have the following lemma.
\begin{lemma}\label{lem:Vassiliev}
The coefficient $c_{p_1,p_2,\dots,p_n}(\nabla_n(L))$ is a Vassiliev invariant
of type $\sum_{i=1}^{n}p_i+1$.
So $C_{p}(\nabla_n(L))$ is also a Vassiliev invariant of type $p+1$.
\end{lemma}
\begin{proof}
Let $\Ell_+$ and $\Ell_-$ be the link diagram as shown below where they are
the same outside this figure.
(This figure has already appeared in R.~Hartley's paper
\cite[Proof of (4.2)]{Hartley:Conway_polynomial}.)
We also let $\Ell^1$ be the singular link diagram which is the
same as $\Ell_+$ and $\Ell_-$ outside the figure.
\vskip 5mm
\begin{gather*}
    \Ell_+:\quad\EllPlus\quad,\qquad
    \Ell_-:\quad\EllMinus\quad,\qquad
    \Ell^1:\quad\EllDot
    \\[5mm]
    \mbox{Figure 2}
\end{gather*}
We assume that $\varphi$ sends $x_1$ and $x_3$ to $t_k$ and $x_2$ and $x_4$
to $t_l$ ($k$ and $l$ may be the same).
Then the Jacobians $M(\Ell_+)$ and $M(\Ell_-)$ are
\begin{equation*}
  M(\Ell_+)
  =
  \left(
    \begin{array}{cccc|c}
      1   &  0  & -1  &  0  & O \\
    1-t_l & t_k &  0  & -1  & O \\
    \hline
     m_1  & m_2 & m_3 & m_4 & M
    \end{array}  
  \right)
\end{equation*}
and
\begin{equation*}
  M(\Ell_-)
  =
  \left(
    \begin{array}{cccc|c}
      1   & t_k-1  & -t_l  &  0  & O \\
      0   &   1    &   0   & -1  & O \\
     \hline
     m_1  & m_2 & m_3 & m_4 & M
    \end{array}
  \right)
\end{equation*}
for some column vectors $m_1,m_2,m_3,m_4$ and a matrix $M$,
where $O$ is a zero vector of suitable size. 
We choose $i,j>4$ to calculate $D^{(ij)}(\Ell_{\pm})$.
Putting
$\tilde D^{(ij)}(\Ell_{\pm})=
D^{(ij)}(\Ell_{\pm})
\prod_{g=1}^{n}{t_g}^{(\kappa_g(\Ell_{\pm})-\nu_g(\Ell_{\pm}))/2}$
and
$\tilde D^{(ij)}(\Ell^1)
=\tilde D^{(ij)}(\Ell_+)-\tilde D^{(ij)}(\Ell_-)$,
we have
\begin{align*}
  \tilde D^{(ij)}(\Ell^1)
  &=
  \prod_{g=1}^{n}{t_g}^{\varepsilon_g}
  \{
  {t_k}^{-1/2}
  \left\Vert
    \begin{array}{cccc|c}
      1  &  0  &  -1  &  0  & O \\
      1  & t_k & -t_l & -1  & O \\
    \hline
    m'_1 & m'_2& m'_3 & m'_4& M'
    \end{array}
  \right\Vert
  \\
  &\phantom{=\prod_{g=1}^{n}{t_g}^{\varepsilon_g}}
  -
  {t_l}^{-1/2}
  \left\Vert
    \begin{array}{cccc|c}
      1   & t_k & -t_l &  -1 & O \\
      0   &  1  &   0  &  -1 & O \\
    \hline
    m'_1  & m'_2& m'_3 & m'_4& M'
    \end{array}
  \right\Vert
  \}
  \\
  &=
  \prod_{g=1}^{n}{t_g}^{\varepsilon_g}
  \{
  \left\Vert
    \begin{array}{cccc|c}
    {t_k}^{-1/2}  &  0  &  -{t_k}^{-1/2}  &  0  & O \\
        1         & t_k &      -t_l       & -1  & O \\
    \hline
      m'_1        & m'_2&       m'_3      & m'_4& M'
    \end{array}
  \right\Vert
  \\
  &\phantom{=\prod_{g=1}^{n}{t_g}^{\varepsilon_g}}
  -
  \left\Vert
    \begin{array}{cccc|c}
      1   &     t_k       & -t_l &      -1       & O \\
      0   & {t_l}^{-1/2}  &   0  & -{t_l}^{-1/2} & O \\
    \hline
     m'_1 &     m'_2      &  m'_3&      m'_4     & M'
    \end{array}
  \right\Vert
  \}
  \\
  &=
  \prod_{g=1}^{n}{t_g}^{\varepsilon_g}
  \left\Vert
    \begin{array}{cccc|c}
      {t_k}^{-1/2} & {t_l}^{-1/2} & -{t_k}^{-1/2} & -{t_l}^{-1/2} & O \\
           1       &     t_k      &     -t_l      &      -1       & O \\
    \hline
          m'_1     &     m'_2     &      m'_3     &      m'_4     & M'
    \end{array}
  \right\Vert
\end{align*}
for some half integers $\varepsilon_1,\varepsilon_2,\dots,\varepsilon_n$,
some column vectors $m'_1,m'_2,m'_3,\linebreak[0]m'_4$, and some matrix $M'$.
(Note that
$\kappa_k(\Ell_+)=\kappa_k(\Ell_-)+1,
 \kappa_l(\Ell_+)=\kappa_l(\Ell_-)-1,
 \nu_k(\Ell_+)=\nu_k(\Ell_-)+2,
 \nu_l(\Ell_+)=\nu_l(\Ell_-)-2$.)
Then the Conway potential function of the labelled singular link $L^1$
presented by $\Ell^1$ is given by
\begin{equation*}
\nabla_n(L^1;t_1,t_2,\dots,t_n)
=\frac{(-1)^{i+j}\tilde D^{(ij)}(\Ell^1)}{\varphi(w_i)(\varphi(x_j)-1)}.
\end{equation*}
\par
Similarly we see that
\begin{equation*}
\nabla_n(L^d;t_1,t_2,\dots,t_n)
=\frac{(-1)^{i+j}\tilde D^{(ij)}(\Ell^d)}{\varphi(w_i)(\varphi(x_j)-1)}.
\end{equation*}
Here $\Ell^d$ is a singular link diagram with $d$ double points presenting
$L^d$ and $\tilde D^{(ij)}(\Ell^d)$ is given as follows.
We arrange $\Ell^d$ so that four arcs adjacent to each double points are
different after inserting kinks if necessary.
Then $\tilde D^{(ij)}(\Ell^d)$ is of the form
\begin{equation*}
\prod_{g=1}^{n}{t_g}^{\varepsilon_g}
\left\Vert
  \begin{array}{ccccc|c}
    T_1  &    O   &   O    & \cdots &    O   &   O   \\
   \ast  &   T_2  &   O    & \cdots &    O   &   O   \\
   \ast  &  \ast  & \ddots & \ddots & \vdots &\vdots \\
  \vdots & \vdots & \ddots & \ddots &    O   &   O   \\
   \ast  &  \ast  & \cdots &  \ast  &   T_d  &   O   \\
   \hline
   \ast  &  \ast  & \cdots &  \ast  &  \ast  & \ast  \\
  \end{array}
\right\Vert,
\end{equation*}
where $T_e$ is a $2\times4$ matrix of the form
\begin{equation*}
\begin{pmatrix}
  {t_{k(e)}}^{-1/2}&{t_{l(e)}}^{-1/2}&-{t_{k(e)}}^{-1/2}&-{t_{l(e)}}^{-1/2} \\
            1      &    t_{k(e)}     &    -t_{l(e)}     &      -1         
\end{pmatrix}.
\end{equation*}
(The author does not know whether the matrix above can be derived from
Fox's free differential calculus applied to $\Ell^d$.)
It is not hard to see that the total degree of $\tilde D^{(ij)}(\Ell^d)$ is
at least $d$ putting $t_k=\exp(h_k)$.
Therefore the total degree of $\nabla_n(L^d)$ is at least $d-1$.
This shows that $c_{p_1,p_2,\dots,p_n}(\nabla_n(L^d))$ vanishes if
$\sum_{k=1}^{n} p_k<d-1$ and so $c_{p_1,p_2,\dots,p_n}(\nabla_n(L))$ is a
Vassiliev invariant of type $\sum_{k=1}^{n}p_k+1$, completing the proof.
\end{proof}
\section{A weight system}
In this section we use J.~Murakami's relations
\cite[p.~126, (1)--(6)]{Jun:multi-variable_Alexander}
to define a weight system $W_n$.
\begin{definition}\label{def:weight}
For a chord diagram $D$ with $d$ chords, we put
$W_{n}(D)\linebreak[0]=C_{d-1}(\nabla_n(D))$ and extend it linearly to a
map from $\D(\stackrel{\mu}{\SplitUnion}S^1)$.
\end{definition}
Since $C_{d-1}$ is a Vassiliev invariant of type $d$ and $D$ has $d$ chords
$W_{n}(D)$ does not depend on its embedding and so is well defined
as a map from $\A(\stackrel{\mu}{\SplitUnion}S^1)$ to the set of homogeneous
polynomials in $h_1,h_2,\dots,\linebreak[0]h_n$
(if $\mu=1$ it also contains a term of the form $c{h_1}^{-1}$).
For a proof that $W_{n}$ satisfies the 4-term relation and the framing
independent relation, see for example
\cite{Birman/Lin:Vassiliev,Birman:new_point}.
\par
Now we will characterise $W_n$.
\begin{proposition}\label{prop:axioms}
$W_{n}$ satisfies the following formulae.
\begin{equation}\label{one_chord}
W_{n}\left(\,\TwoParallelArcsUpUpWithOneChord{$i$}{$i$}\,\right)
=h_iW_{n}\left(\,\TwoCrossedArcsUpUp{$i$}{$i$}\,\right),
\end{equation}
\medskip
\begin{multline}\label{Jun's_relation_for_chord_diagram}
  4h_j
  \{
    W_{n}\left(\,\AmidaBC{$i$}{$j$}{$k$}\,\right)
   -W_{n}\left(\,\AmidaCA{$i$}{$j$}{$k$}\,\right)
  \}
  \\[5mm]
  +2(h_i-h_k)
  \{
    W_{n}\left(\,\AmidaAB{$i$}{$j$}{$k$}\,\right)
   +W_{n}\left(\,\AmidaBA{$i$}{$j$}{$k$}\,\right)
  \}
  \\[5mm]
  +(h_k-h_i)(h_ih_k+{h_j}^2)
  W_{n}\left(\,\AmidaI{$i$}{$j$}{$k$}\,\right)
  =0,
\end{multline}
\medskip
\begin{equation}\label{one_arc_one_circle}
W_{n}\left(\,\OneArcUpAndOneCircleWithOneChord{$i$}{$j$}\,\right)
=h_iW_{n}\left(\,\OneArcUp{$i$}\,\right),
\end{equation}
\medskip
\begin{equation}\label{one_circle}
W_{n}\left(\,\OneCircle{$i$}\,\right)
={h_i}^{-1},
\end{equation}
\medskip
\begin{equation}\label{split}
W_{n}
\left(\mbox{\rm any nonempty chord diagram}\,\SplitUnion\,\OneCircle{}\right)
=0.
\end{equation}
Here $i,j,k$ indicate the labels attached to the arcs near by.
Note that some of the labels $i,j,k$ may be equal.  
Note also that the crossing in the right hand side of (\ref{one_chord})
is {\em not} a double point.
It only indicates the connectivity. 
\end{proposition}
\begin{proof}
We assume that the chord diagrams appearing in the lemma have $d$ double
points outside the regions described in the pictures.
\par
(Proof of (\ref{one_chord}))
From the well-known relation for the potential function
(Conway's first identity, which is
the first relation of \cite[p126]{Jun:multi-variable_Alexander})
\begin{equation*}
\nabla_n\left(\,\PositiveCrossUpUp{$i$}{$i$}\,\right)
-
\nabla_n\left(\,\NegativeCrossUpUp{$i$}{$i$}\,\right)
=
2\sinh(h_i/2)\nabla_n\left(\,\TwoParallelArcsUpUp{}{$i$}{$i$}\,\right),
\end{equation*}
we have
\begin{equation*}
\nabla_n\left(\,\PositiveFullTwistUpUp{$i$}{$i$}\,\right)
-
\nabla_n\left(\,\TwoParallelArcsUpUp{}{$i$}{$i$}\,\right)
=
2\sinh(h_i/2)\nabla_n\left(\,\PositiveCrossUpUp{$i$}{$i$}\,\right).
\medskip
\end{equation*}
Note that this also holds for singular links and recall that we are assuming
that there are $d$ double points outside the region appearing in the equality
above.
Taking the total degree $d$ part, we have
\begin{align}\label{eq:crossing_change}
  C_{d}(\nabla_n\left(\,\PositiveFullTwistUpUp{$i$}{$i$}\,\right))
  -
  C_{d}(\nabla_n\left(\,\TwoParallelArcsUpUp{}{$i$}{$i$}\,\right))
  &=
  2C_{1}(\sinh(h_i/2))C_{d-1}
  (\nabla_n\left(\,\PositiveCrossUpUp{$i$}{$i$}\,\right))
  \\[5mm]
  &=
  h_iC_{d-1}(\nabla_n\left(\,\PositiveCrossUpUp{$i$}{$i$}\,\right))
  \notag
\end{align}
since 
$C_{e}(\nabla_n\left(\OneArcUp{$i$}\right)$\vspace{5mm}
vanishes if $e<d-1$ from Lemma~\ref{lem:Vassiliev}.
By the way, we have from the definition
\begin{align*}
  W_{n}\left(\,\TwoParallelArcsUpUpWithOneChord{$i$}{$i$}\,\right)
  &=
  C_{d}(\nabla_n\left(\,\TwoBraidCP{$i$}{$i$}\,\right))
  \\[5mm]
  &=
  C_{d}(\nabla_n\left(\,\PositiveFullTwistUpUp{$i$}{$i$}\,\right))
  -
  C_{d}(\nabla_n\left(\,\TwoParallelArcsUpUp{}{$i$}{$i$}\,\right))
\end{align*}
and
\begin{equation*}
  W_{n}\left(\,\TwoCrossedArcsUpUp{$i$}{$i$}\,\right)
  =
  C_{d-1}(\nabla_n\left(\,\PositiveCrossUpUp{$i$}{$i$}\,\right)).
\medskip
\end{equation*}
Therefore the required formula follows from (\ref{eq:crossing_change}).
\par
(Proof of (\ref{one_arc_one_circle}))
From the fourth relation of \cite[p.126]{Jun:multi-variable_Alexander} we have
\begin{equation*}
  \nabla_n\left(\,\PositiveClasp{$i$}{$j$}\,\,\right)
  =
  2\sinh(h_i/2)\nabla_n\left(\OneArcUp{$i$}\right).
\medskip
\end{equation*}
Taking the total degree $d$ part of both hand sides, we have
\begin{align}\label{clasp}
  C_{d}(\nabla_n\left(\,\PositiveClasp{$i$}{$j$}\,\,\right))
  &=
  2C_{1}(\sinh(h_i/2))C_{d-1}(\nabla_n\left(\OneArcUp{$i$}\right))
  \\[5mm]
  &=
  h_iC_{d-1}(\nabla_n\left(\OneArcUp{$i$}\right)).
  \notag
\end{align}
So we have
\begin{align*}
  W_{n}\left(\,\OneArcUpAndOneCircleWithOneChord{$i$}{$j$}\,\right)
  &=
  C_{d}(\nabla_n\left(\,\OneArcUpAndOneCircleWithOneChord{$i$}{$j$}\,\right))
  \\[5mm]
  &=
  C_{d}(\nabla_n\left(\,\PositiveClaspWithDoublePoint{$i$}{$j$}\,\,\right))
  \\[5mm]
  &=
  C_{d}(\nabla_n\left(\,\PositiveClasp{$i$}{$j$}\,\,\right))
  -C_{d}(\nabla_n\left(\,\OneArcUpAndOneCircleWithoutChord{$i$}{$j$}\,\right))
  \\[5mm]
  &=
  h_iC_{d-1}(\nabla_n\left(\OneArcUp{$i$}\right))
  \\[5mm]
  &=
  h_iW_{n}\left(\OneArcUp{$i$}\right).
\end{align*}
Here we use (\ref{clasp}) and the fact that $\nabla_n$ vanishes for a split
link in the fourth equality.
\par
(Proof of (\ref{one_circle}) and (\ref{split}))
Since $\nabla_n(O_{i})=1/(2\sinh(h_i/2))$ (which is the 
fifth relation of \cite[p.~126]{Jun:multi-variable_Alexander}),
$W_{n}(O_{i})=C_{-1}(1/(2\sinh(h_i/2)))={h_i}^{-1}$
and we have (\ref{one_circle}).
Here $O_{i}$ is the trivial knot with label $i$.
The relation (\ref{split}) follows from the fact that $\nabla_n$ vanishes for
split links (the sixth relation of
\cite[p.~126]{Jun:multi-variable_Alexander}).
\par
(Proof of (\ref{Jun's_relation_for_chord_diagram}))
We use J.~Murakami's third relation 
\cite[p.126]{Jun:multi-variable_Alexander}:
\begin{align*}
  &4\cosh(h_i/2)\sinh(h_j/2)
    \nabla_n\left(\,\ThreeBraidPTwoPOnePOnePTwo{$i$}{$j$}{$k$}\,\right)
  \\[5mm]
  -&4\cosh(h_k/2)\sinh(h_j/2)
    \nabla_n\left(\,\ThreeBraidPOnePTwoPTwoPOne{$i$}{$j$}{$k$}\,\right)
  \\[5mm]
  -&2\sinh((-h_i+h_k)/2)
    \{
      \nabla_n\left(\,\ThreeBraidPTwoPTwoPOnePOne{$i$}{$j$}{$k$}\,\right)
     +\nabla_n\left(\,\ThreeBraidPOnePOnePTwoPTwo{$i$}{$j$}{$k$}\,\right)
    \}
  \\[5mm]
  +&4\cosh(h_k/2)\sinh((-h_i+h_j+h_k)/2)
    \nabla_n\left(\,\ThreeBraidPOnePOne{$i$}{$j$}{$k$}\,\right)
  \\[5mm]
  -&4\cosh(h_i/2)\sinh(( h_i+h_j-h_k)/2)
    \nabla_n\left(\,\ThreeBraidPTwoPTwo{$i$}{$j$}{$k$}\,\right)
  \\[5mm]
  -&2\sinh(-h_i+h_k)
    \nabla_n\left(\,\ThreeParallelArcsUpUpUp{}{$i$}{$j$}{$k$}\,\right)
  =0.
\medskip
\end{align*}
We take the total degree $d+2$ part.
(Recall that we assume that there are $d$ double points outside.)
Since $C_{e}(L^d)=0$ for $e<d-1$ and any singular link $L^d$ with $d$ double
points, we have
\begin{align*}
  &4C_{1}(\cosh(h_i/2)\sinh(h_j/2))
    C_{d+1}
      (\nabla_n\left(\,\ThreeBraidPTwoPOnePOnePTwo{$i$}{$j$}{$k$}\,\right))
  \\[5mm]
  &-4C_{1}(\cosh(h_k/2)\sinh(h_j/2))
     C_{d+1}
       (\nabla_n\left(\,\ThreeBraidPOnePTwoPTwoPOne{$i$}{$j$}{$k$}\,\right))
  \\[5mm]
  &-2C_{1}(\sinh((-h_i+h_k)/2))
    \{
      C_{d+1}
        (\nabla_n\left(\,\ThreeBraidPTwoPTwoPOnePOne{$i$}{$j$}{$k$}\,\right))
     +C_{d+1}
        (\nabla_n\left(\,\ThreeBraidPOnePOnePTwoPTwo{$i$}{$j$}{$k$}\,\right))
    \}
  \\[5mm]
  &+4C_{1}(\cosh(h_k/2)\sinh((-h_i+h_j+h_k)/2))
     C_{d+1}(\nabla_n\left(\,\ThreeBraidPOnePOne{$i$}{$j$}{$k$}\,\right))
  \\[5mm]
  &-4C_{1}(\cosh(h_i/2)\sinh(( h_i+h_j-h_k)/2))
     C_{d+1}(\nabla_n\left(\,\ThreeBraidPTwoPTwo{$i$}{$j$}{$k$}\,\right))
  \\[5mm]
  &-2C_{1}(\sinh(-h_i+h_k))
     C_{d+1}
       (\nabla_n\left(\,\ThreeParallelArcsUpUpUp{}{$i$}{$j$}{$k$}\,\right))
  \\[5mm]
  &+4C_{3}(\cosh(h_i/2)\sinh(h_j/2))
     C_{d-1}
       (\nabla_n\left(\,\ThreeBraidPTwoPOnePOnePTwo{$i$}{$j$}{$k$}\,\right))
  \\[5mm]
  &-4C_{3}(\cosh(h_k/2)\sinh(h_j/2))
     C_{d-1}
       (\nabla_n\left(\,\ThreeBraidPOnePTwoPTwoPOne{$i$}{$j$}{$k$}\,\right))
  \\[5mm]\pagebreak[0]
  &-2C_{3}(\sinh((-h_i+h_k)/2))
    \{
      C_{d-1}
        (\nabla_n\left(\,\ThreeBraidPTwoPTwoPOnePOne{$i$}{$j$}{$k$}\,\right))
     +C_{d-1}
        (\nabla_n\left(\,\ThreeBraidPOnePOnePTwoPTwo{$i$}{$j$}{$k$}\,\right))
    \}
  \\[5mm]
  &+4C_{3}(\cosh(h_k/2)\sinh((-h_i+h_j+h_k)/2))
     C_{d-1}(\nabla_n\left(\,\ThreeBraidPOnePOne{$i$}{$j$}{$k$}\,\right))
  \\[5mm]
  &-4C_{3}(\cosh(h_i/2)\sinh(( h_i+h_j-h_k)/2))
     C_{d-1}(\nabla_n\left(\,\ThreeBraidPTwoPTwo{$i$}{$j$}{$k$}\,\right))
  \\[5mm]
  &-2C_{3}(\sinh(-h_i+h_k))
     C_{d-1}
       (\nabla_n\left(\,\ThreeParallelArcsUpUpUp{}{$i$}{$j$}{$k$}\,\right))
  =0.
\end{align*}
\vspace{5mm}
\par
Now since $C_{d-1}$ is a type $d$ invariant,\vspace{5mm}
all its values in the equality above are the same and equal to
$C_{d-1}(\nabla_n\left(\,\ThreeParallelArcsUpUpUp{}{$i$}{$j$}{$k$}\,\right))$.
So we have
\vspace{5mm}
\begin{align}\label{degree_d_of_Jun's_relation}
  &2h_j
    \{
      C_{d+1}
        (\nabla_n\left(\,\ThreeBraidPTwoPOnePOnePTwo{$i$}{$j$}{$k$}\,\right))
     -C_{d+1}
        (\nabla_n\left(\,\ThreeBraidPOnePTwoPTwoPOne{$i$}{$j$}{$k$}\,\right))
    \}
  \\[5mm]
  &-(-h_i+h_k)
    \{
      C_{d+1}
        (\nabla_n\left(\,\ThreeBraidPTwoPTwoPOnePOne{$i$}{$j$}{$k$}\,\right))
     +C_{d+1}
        (\nabla_n\left(\,\ThreeBraidPOnePOnePTwoPTwo{$i$}{$j$}{$k$}\,\right))
    \}
  \notag
  \\[5mm]
  &+2(-h_i+h_j+h_k)
    C_{d+1}(\nabla_n\left(\,\ThreeBraidPOnePOne{$i$}{$j$}{$k$}\,\right))
  \notag
  \\[5mm]
  &-2(h_i+h_j-h_k)
    C_{d+1}(\nabla_n\left(\,\ThreeBraidPTwoPTwo{$i$}{$j$}{$k$}\,\right))
  \notag
  \\[5mm]
  &-2(-h_i+h_k)
    C_{d+1}(\nabla_n\left(\,\ThreeParallelArcsUpUpUp{}{$i$}{$j$}{$k$}\,\right))
  \notag
  \\[5mm]
  &+\frac{1}{2}(h_k-h_i)(h_ih_k+h_j^2)
    C_{d-1}(\nabla_n\left(\,\ThreeParallelArcsUpUpUp{}{$i$}{$j$}{$k$}\,\right))
  =0.
  \notag
\end{align}
Now we have from the definition of $W_n$
\begin{align*}
  &2h_j
  \{
    W_n\left(\,\AmidaBC{$i$}{$j$}{$k$}\,\right)
    -
    W_n\left(\,\AmidaCA{$i$}{$j$}{$k$}\,\right)
  \}
  \\[5mm]
  &+(h_i-h_k)
  \{
    W_n\left(\,\AmidaAB{$i$}{$j$}{$k$}\,\right)
    +
    W_n\left(\,\AmidaBA{$i$}{$j$}{$k$}\,\right)
  \}
  \\[5mm]
  &+\frac{1}{2}(h_k-h_i)(h_ih_k+{h_j}^2)
  W_{n}\left(\,\AmidaI{$i$}{$j$}{$k$}\,\right)
  \\[5mm]
  =
  &2h_j
  \{
   C_{d+1}
     (\nabla_n\left(\,\ThreeBraidPTwoPOneCOneCTwo{$i$}{$j$}{$k$}\,\right))
   -
   C_{d+1}
     (\nabla_n\left(\,\ThreeBraidCOneCTwoPTwoPOne{$i$}{$j$}{$k$}\,\right))
  \}
  \\[5mm]
  &+(h_i-h_k)
  \{
    C_{d+1}
      (\nabla_n\left(\,\ThreeBraidPTwoCTwoPOneCOne{$i$}{$j$}{$k$}\,\right))
    +
    C_{d+1}
      (\nabla_n\left(\,\ThreeBraidPOneCOnePTwoCTwo{$i$}{$j$}{$k$}\,\right))
  \}
  \\[5mm]
  &+\frac{1}{2}(h_k-h_i)(h_ih_k+{h_j}^2)
    C_{d+1}(\nabla_n\left(\,\AmidaI{$i$}{$j$}{$k$}\,\right)),
  \\[5mm]
\end{align*}
which vanishes from (\ref{degree_d_of_Jun's_relation}) and the proof is
complete.
\end{proof}
\begin{remark}
In the proof above we did not use Conway's second identity (the second
relation of \cite[p.~126]{Jun:multi-variable_Alexander}.
The author does not know whether it is necessary in J.~Murakami's axioms for
the multivariable Alexander polynomial.
In our case we have the following corollary which corresponds to Conway's
second identity.
\end{remark}
\begin{corollary}\label{cor:two_chords}
We have
\begin{align*}
  W_{n}\left(\,\TwoParallelArcsUpUpWithTwoParallelChords{$i$}{$j$}\,\right)
  &=
  \left(\frac{i+j}{2}\right)^2
  W_{n}\left(\,\TwoParallelArcsUpUp{}{$i$}{$j$}\,\right)
\end{align*}
and
\begin{align*}
  W_{n}\left(\,\TwoParallelArcsUpUpWithTwoCrossedChords{$i$}{$j$}\,\right)
  &=
  \left(\frac{i-j}{2}\right)^2
  W_{n}\left(\,\TwoParallelArcsUpUp{}{$i$}{$j$}\,\right).
\end{align*}
\end{corollary}
\begin{proof}
We put $k=j$ in the relation (\ref{Jun's_relation_for_chord_diagram}) and
connect these two arcs as follows.
\begin{align*} 
  &4h_j
  \{W_{n}\left(\,\AmidaBCWithCurlTop{$i$}{$j$}\,\right)
    -
    W_{n}\left(\,\AmidaCAWithCurlTop{$i$}{$j$}\,\right)
  \}
  \\[5mm]
  &+2(h_i-h_j)
  \{
    W_{n}\left(\,\AmidaABWithCurlTop{$i$}{$j$}\,\right)
    +
    W_{n}\left(\,\AmidaBAWithCurlTop{$i$}{$j$}\,\right)
  \}
  \\[5mm]
  &+h_j(h_j-h_i)(h_i+h_j)
  W_{n}\left(\,\AmidaIWithCurlTop{$i$}{$j$}\,\right)
  =0.
\end{align*}
\vspace{5mm}
Then applying the relation (\ref{one_chord}), we have
\begin{align*}
  &4h_j
  \{
    h_jW_{n}
      \left(\,
        \TwoArcsUpUpAndOneCircleWithOneChordToCircle{$i$}{$j$}{$j$}
      \,\right)
    -
    W_{n}
      \left(\,
        \TwoParallelArcsUpUpWithTwoParallelChords{$i$}{$j$}
      \,\right)
  \}
  \\[5mm]
  &+2h_j(h_i-h_j)
  \{
    W_{n}
      \left(\,
        \TwoArcsUpUpAndOneCircleWithOneChordToCircle{$i$}{$j$}{$j$}
      \,\right)
    +
    W_{n}
      \left(\,
        \TwoArcsUpUpAndOneCircleWithOneChordToArcs{$i$}{$j$}{$j$}
      \,\right)
  \}
  \\[5mm]
  &+h_j(h_j-h_i)(h_i+h_j)
  W_{n}
    \left(\,
      \TwoParallelArcsUpUp{}{$i$}{$j$}
    \,\right)
  =0.
  \\[5mm]
\end{align*}
Now using the relations (\ref{one_arc_one_circle}) and (\ref{split}), we have
\begin{align*}
  \{4h_j^2h_i+h_j(h_j-h_i)(h_i+h_j)+2h_ih_j(h_i-h_j)\}
  &W_{n}\left(\,\TwoParallelArcsUpUp{}{$i$}{$j$}\,\right)
  \\[5mm]
  &-4h_j
  W_{n}\left(\,\TwoParallelArcsUpUpWithTwoParallelChords{$i$}{$j$}\,\right)
  =0.
\end{align*}
\vspace{5mm}
So the required formula follows.
\par
Similarly the following connection shows the second formula.
\begin{align*} 
  &4h_j
  \{
    W_{n}\left(\,\AmidaBCWithCurlBottom{$i$}{$j$}\,\right)
    -
    W_{n}\left(\,\AmidaCAWithCurlBottom{$i$}{$j$}\,\right)
  \}
\\[5mm]
  &+2(h_i-h_j)
  \{
    W_{n}\left(\,\AmidaABWithCurlBottom{$i$}{$j$}\,\right)
    +
    W_{n}\left(\,\AmidaBAWithCurlBottom{$i$}{$j$}\,\right)
  \}
  \\[5mm]
  &+h_j(h_j-h_i)(h_i+h_j)
  W_{n}\left(\,\AmidaIWithCurlBottom{$i$}{$j$}\,\right)
  =0.
\end{align*}
\end{proof}
\par
J.~Murakami proved that his six relations are sufficient to calculate the
multivariable Alexander polynomial.
Now our main result is
\begin{theorem}
$W_n$ can be calculated recursively by using axioms (\ref{one_chord}) --
(\ref{split}).
\end{theorem}
\begin{proof}
We proceed by induction on the number of circles in the support of a chord
diagram.
If there is only one circle, then we use (\ref{one_chord}) to change the
chord diagram into a diagram without chords.
Then we apply (\ref{one_circle}) or (\ref{split}) to evaluate the diagram.
\par
Suppose that we are given a chord diagram $D$ with support
$E_1,E_2,\dots,E_{\mu}$ ($\mu>1$).
We first look at $E_1$.
If $E_1$ contains no end point of a chord, then $W_n(D)=0$ from (\ref{split}).
If $E_1$ contains one end point, then from (\ref{one_arc_one_circle})
we can reduce the number of circles.
If $E_1$ contains more than one end point, we use
(\ref{Jun's_relation_for_chord_diagram}).
We assume that $E_1$ is labelled $i$ in
(\ref{Jun's_relation_for_chord_diagram}).
The second term there contains two end points and the others contain one or
less.
So we can reduce the number of end points.
Repeating this process we have chord diagrams with one or no end point, which
can be calculated as described above.
\par
So the proof is complete.
\end{proof}
\section{Problems and a conjecture}
In this section we discuss open problems and state a conjecture.
Our first problem is
\begin{problem}
Show that $W_{n}$ is well-defined without using the Conway potential function.
\end{problem}
Note that the relation (\ref{Jun's_relation_for_chord_diagram}) implies the
4-term relation as follows.
We can write (\ref{Jun's_relation_for_chord_diagram}) as
\begin{equation*}
4h_j(BC-CA)+2(h_i-h_k)(AB+BA)+(h_k-h_i)(h_ih_k+{h_j}^2)I=0
\end{equation*}
with $I=\AmidaISmall{$i$}{$j$}{$k$}$\,, $A=\AmidaA{$i$}{$j$}{$k$}$\,,
$B=\AmidaB{$i$}{$j$}{$k$}$\,, and $C=\AmidaC{$i$}{$j$}{$k$}$ and we composite
them downward.
If we exchange the labels $i$ and $k$, then $A$ and $B$ are also exchanged
and so we have
\begin{equation*}
4h_j(AC-CB)+2(h_k-h_i)(AB+BA)+(h_i-h_k)(h_ih_k+{h_j}^2)I=0.
\end{equation*}
Adding the two equalities above and divide by $4h_j$, we have
$AC-CA=CB-BC$, which is the 4-term relation.
\par
It is easily seen that (\ref{one_chord}) and (\ref{split}) imply the
framing independence relation.
Therefore the well-definedness of $W_{n}$ as a map from
$\D(\SplitUnion S^1)$ implies that it factors through $\A(\SplitUnion S^1)$,
which proves that $W_{n}$ defines a link invariant via the Kontsevich
integral!
\par
The next problem is
\begin{problem}
Can we alter the coefficients appeared in Proposition~\ref{prop:axioms}?
\end{problem}
For example, let us replace (\ref{Jun's_relation_for_chord_diagram}) with
\begin{equation}\label{5-term}
x(h_i,h_j,h_k)(BC-CA)+y(h_i,h_j,h_k)(AB+BA)+z(h_i,h_j,h_k)I=0,
\end{equation}
where $x(h_i,h_j,h_k),y(h_i,h_j,h_k)$ and $z(h_i,h_j,h_k)$ are functions of
$h_i,h_j,$ and $h_k$.
If $x(h_i,h_j,h_k)$ is (nonzero and) symmetric with respect to $h_i$ and
$h_k$, and $y(h_i,h_j,h_k)$ and $z(h_i,h_j,h_k)$ are antisymmetric with
respect to $h_i$ and $h_k$, then the relation (\ref{5-term}) above implies
the 4-term relation.
So if we could prove that $\D(\SplitUnion S^1)$ modulo (\ref{5-term}),
(\ref{one_chord}), (\ref{one_arc_one_circle}), (\ref{one_circle}), and
(\ref{split}) is nontrivial, we might have another labelled link invariant.
Our $W_n$ is a map to homogeneous polynomials but the author does not know
whether the homogeneity is necessary or not.
\par
As D.~Bar-Natan and S.~Garoufalidis pointed out in
\cite{BarNatan/Garoufalidis:MMR_conjecture}, $W_{1}$ is a canonical
weight system in the sense that $W_{1}\circ Z$ coincides with
$\nabla_1$.
(Note that their definition of the weight system is $h_1W_{1}$ in our
notation.)
Our conjecture is
\begin{conjecture}
$W_{n}$ is canonical for every $n$, i.e., $W_{n}\circ Z=\nabla_n$. 
\end{conjecture}
To prove this it is sufficient to prove that $W_n\circ Z$ satisfies the six
axioms in \cite[p.~126]{Jun:multi-variable_Alexander}.
The author cannot prove it since they involve Drinfel'd's associator
\cite{Drinfeld:quasi-Hopf_algebras,%
Drinfeld:Gal(Q/Q),%
Le/Murakami:Kauffman_polynomial,%
Le/Murakami:q-tangle,%
BarNatan:non-associative_tangles,%
Kassel:quantum_groups}.
Note that there is no direct proof (using Drinfel'd's associator) of the
equality $W_1\circ Z(O_1)=1/(2\sinh(h_1))$.
The proof in \cite[Example 2.7]{BarNatan/Garoufalidis:MMR_conjecture}
depends on T.Q.T.~Le and J.~Murakami's result on the canonical weight
system for the HOMFLY polynomial
\cite{Le/Murakami:HOMFLY_polynomial} which depends on the skein relation
not on Drinfel'd's associator
(in fact they use the skein relation to prove some interesting formulae
about coefficients of Drinfel'd's associator which involve the multiple
zeta functions).
Note also that since $W_{1}=W_{n}\big\vert_{h_1=h_2=\cdots=h_n}$, we have
$(W_{n}\circ Z)\big\vert_{h_1=h_2=\cdots=h_n}
=\nabla_n\big\vert_{h_1=h_2=\cdots=h_n}$. 
\bibliography{journal,hitoshi}
\bibliographystyle{amsplain}
\end{document}